\documentclass[12pt]{iopart}
\usepackage{iopams,setstack}

\usepackage[latin1]{inputenc}

\providecommand{\U}[1]{\protect\rule{.1in}{.1in}}

\providecommand{\U}[1]{\protect\rule{.1in}{.1in}}
\providecommand{\U}[1]{\protect\rule{.1in}{.1in}} \makeatletter
\providecommand{\U}[1]{\protect\rule{.1in}{.1in}}

\newtheorem{theorem}{Theorem}

\newtheorem{corollary}{Corollary}

\newtheorem{definition}{Definition}
\newtheorem{example}{Example}

\newtheorem{lemma}{Lemma}

\newtheorem{proposition}{Proposition}
\newtheorem{remark}{Remark}

\newenvironment{proof}[1][Proof]{\noindent\textbf{#1.} }{\ \rule{0.5em}{0.5em}}
\def\dfrac{\frac}

\begin{document}

\title[Differential invariants of generic PMAs]
 {Differential invariants of generic parabolic Monge-Amp\`ere
equations} 

\author{D. Catalano Ferraioli$^1$ and A. M. Vinogradov$^2$}

\address{$^1$ Departamento de Matem\' atica, Universidade Federal da Bahia, Av. Ademar de Barros s/n, 40170-110 Salvador, BA, Brasil}
\address{$^2$ Dipartimento di Matematica e Informatica Universit\`a degli Studi di Salerno, via Ponte Don Melillo, 84084 Fisciano (SA), Italia}

\ead{diego.catalano@ufba.br and vinograd@unisa.it}

\begin{abstract}
Some new results on geometry of classical parabolic Monge-Amp\`ere
equations (PMA) are presented. PMAs are either \emph{integrable}, or
\emph{nonintegrable} according to integrability of its
characteristic distribution. All integrable PMAs are locally
equivalent to the equation $u_{xx}=0$. We study nonintegrable PMAs
by associating with each of them a $1$-dimensional distribution on
the corresponding first order jet manifold, called the
\emph{directing distribution}. According to some property of this
distribution, nonintegrable PMAs are subdivided into three
classes, one \emph{generic} and two \emph{special} ones. Generic
PMAs are completely characterized by their directing distributions,
and we study canonical models of the latters, \emph{projective curve 
bundles} (PCB). A PCB is a
$1$-dimensional subbundle of the projectivized cotangent bundle of
a $4$-dimensional manifold. Differential invariants of projective
curves composing such a bundle are used to construct a series of
contact differential invariants for corresponding PMAs. These give
a solution of the equivalence problem for generic PMAs with respect to
contact transformations. The introduced invariants measure in an exact
manner nonlinearity of PMAs.
\end{abstract}

\pacs{02.20.Qs, 02.20.TW, 02.30.Jr}
\ams{34A26, 53A55,  53D10, 53A20, 35K55}
\noindent{\it Keywords\/}: Differential invariants, jet spaces, Monge-Amp\`ere equations, parabolic equations, contact transformations

\section{Introduction}

After the classical ``Application de l'Analyse \`{a} la
G\'{e}om\'{e}trie'' by G.~Monge, Monge-Amp\`ere equations continuously     
attract attention of geometers and physicists. Indeed they are intimately related with numerous parts of modern
differential geometry, like for example Calabi-Yau conjecture, but also have various applications  to physics. For example, these equations appears in cosmological models and describe classes of exact solutions of Einstein's field equations (see for example \cite{NutShef,Pavlov} and the extensive account in \cite{Stephani}). Other applications of these equations can also be found in \cite{Gutierrez,KLR,Nutku}.\\
Last 2-3 decades manifested a new wave of interest in the geometry of
these equations, mostly of elliptic and hyperbolic types.
In \cite{KLR} the reader will find a rather complete account of classical and recent 
results on the geometry of Monge-Amp\`ere equations and also an
extensive bibliography.

In this article we study geometry of classical parabolic
Monge-Ampere equations (PMAs) on the basis of a new approach
sketched in \cite{V}. It is based on the observation that a PMA $\mathcal{E}\subset
J^{2}(\pi), \,\pi$ being a 1-dimensional fiber bundle over a
bidimensional manifold,  is completely characterized by 
a 2-dimensional Lagrangian distribution $\mathcal{D_{\mathcal{E}}}$
on $J^{1}(\pi)$ , called the \emph{characteristic distribution} of
$\mathcal{E}$, and vice versa.
 Such distributions and, accordingly, the
corresponding PMAs, are naturally subdivided into four
classes of \emph{integrable, generic} and two 
\emph{special} types of equations (see \cite{V} and sec.4). 
In this classification, integrable PMAs are those whose characteristic distributions are integrable. 
Since all integrable Lagrangian distributions are locally contact equivalent,
all integrable PMAs are locally contact equivalent to one of them, say,
to the equation $u_{xx}=0$.
On the contrary, nonintegrable Lagrangian
distributions are very diversified and our main goal here is to
describe their multiplicity, i.e., equivalence
classes of PMAs with respect to contact transformations.

With this purpose we associate with a Lagrangian distribution a
\emph{projective curve bundle} (shortly, PCB) over a 4-dimensional
manifold $N$. A PCB over $N$ is a 1-dimensional smooth subbundle
of the ``projectivized'' cotangent bundle $PT^{\ast}(N)$ of $N$.
Under some regularity conditions such a bundle possesses a
canonical contact structure and, as a consequence, a Lagrangian distribution
canonically inscribed in it. There exists a one-to-one
correspondence between generic Lagrangian distributions and
\emph{regular} PCBs. This way the equivalence problem for generic
PMAs is reformulated as the equivalence problem for PCBs (see
\cite{V}). This is the key point of our approach. The fiber of
such a bundle over a point $x\in N$ is a curve $\gamma_{x}$ in the
projective space $PT_{x}^{\ast}(N)$. The curve $\gamma_{x}$ can be
characterized by its scalar differential invariants with respect
to the group of projective transformations. By putting together these
invariants for single curves $\gamma_{x}, \,x\in
N$, one obtains scalar differential invariants for the considered
PCB and, consequently, for the corresponding PMA.

This kind of invariants resolves the equivalence problem for
generic PMAs on the basis of the "principle of $n-$invariants"
(see \cite{AVL,SDI}). We have focused on these invariants because of their
transparent geometrical meaning. However, it should be stressed
that there are other choices, maybe, not equally clear geometrically, but more 
computationally convenient. They will be discussed separately.

Special Lagrangian distributions admit a similar interpretation in
terms of 2-dimensional distributions on 4-dimensional manifolds
supplied with additional structures, called \emph{fringes} (see
\cite{V}). Differential invariants of fringes, also coming from
projective differential geometry, allow to construct basic scalar
differential invariants for special PMAs. They will be discussed
in a separate paper. It is worth mentioning that linear 
 as well as many other concrete PMAs of current interest are very 
 "simple" in comparison with generic ones. Rigorously, this means that
 they belong to the special class.

Our approach is motivated by the theory of solution singularities for
nonlinear PDEs (see \cite{VinSing}) and has been already applied to hyperbolic
and elliptic Monge-Amp\`ere equations (see \cite{MVY,DPV}). Namely, the characteristic 
distribution $\mathcal{D_{\mathcal{E}}}$ of a PMA $\mathcal{E}$
or, equivalently, the associated PCB, is, in fact, the
equation that describe fold type singularities of multivalued
solutions of $\mathcal{E}$. An advantage of this approach
is that it naturally extends to higher order PDEs
and, in particular, allows one to distinguish higher order analogues
of Monge-Amp\`ere equations. These topics will be discussed in a paper 
in preparation by M. Bachtold and the second author.

The paper is organized as follows. The notations and necessary 
generalities concerning jet spaces and Monge-Amp\`ere equations 
are collected in sections 2 and 3,
respectively. In particular, the interpretation of PMAs as
Lagrangian distributions is explained in sec.3.  The central notion of the paper, 
namely, that of the \emph{directing distribution}  
of a Lagrangian distribution, and other basic facts of 
geometry of the Lagrangian distributions are discussed in section 4.
The aforementioned subdivision of
nonitegrable PMAs into generic and special types reflects some
contact properties of this distribution. For completeness
in section 5 we give a short proof of the known fact that
integrable PMAs are locally equivalent each other.
Projective curve bundles, canonical models of directing distributions, are
introduced and studied in section 6. Finally,
basic scalar differential invariants of generic PCBs and hence of
generic PMAs are constructed and discussed in section 7.

Throughout the paper we use the following notations and
conventions:

\begin{itemize}
\item all objects in this paper, e.g., manifolds, mappings,
functions, vector fields, etc, are supposed to be smooth;

\item  $C^{\infty}(M)$ stands for the algebra of smooth functions
on the manifold $M$ and $C^{\infty}(M)$-modules of all vector
fields and differential $k$-forms are denoted by $D(M)$ and
$\Lambda^{k}(M)$, respectively;

\item the evaluation of $X\in D(M)$ (respectively, of
$\alpha\in\Lambda^{k}(M)$) at $p\in M$ is denoted by $\left.
X\right\vert _{p}$ ( respectively, $\left.  \alpha\right\vert _{p}$);

\item  $d_{p}f:T_{p}M\longrightarrow T_{f(p)}N$ stands for the
differential of the map $f:M\rightarrow N$ at $p$;

\item For $X,Y\in D(M)$ and $\alpha\in\Lambda^{k}(M)$ we  use
the short notation $X^{r}(Y)$ and $X^{r}(\alpha)$  for $L_{X}^r(Y)$ 
and $L_{X}^r(\alpha)$  (the $r$-th power of the Lie derivative
$L_{X}$), respectively, assuming that $X^{0}(Y)=Y$ and
$X^{0}(\alpha)=\alpha$.

\item We slightly abuse the language by using the term \emph{distribution} on a
manifold $M$  either for a subbundle $\mathcal{D}$ of the
tangent bundle  $TM$, whose fiber over $p\in M$ is denoted by
$\mathcal{D}(p)$, or for the $C^{\infty}(M)$-module of its sections.
In particular, $X\in\mathcal{D}$ means that
$X_{p}\in\mathcal{D}(p)$, $\forall p\in M$;

\item we write $\mathcal{D}=\left\langle
X_{1},...,X_{2}\right\rangle $ if the distribution $\mathcal{D}$
is generated by vector fields $X_{1},...,X_{r}\in D(M)$;
similarly, $\mathcal{D}=\mathrm{Ann}(\alpha _{1},...,\alpha_{s})$ means
that $\mathcal{D}$ is constituted by vector fields annihilated by
forms $\alpha_{1},...,\alpha_{s}\in\Lambda^{1}(M)$;

\item if $M$ is equipped with a contact distribution $\mathcal{C}$
and $\mathcal{S}\subset\mathcal{C}$ is a subdistribution of
$\mathcal{C}$, then $\mathcal{S}^{\bot}$ denotes the
$\mathcal{C}$-orthogonal complement to $\mathcal{S}$.
\end{itemize}

\section{Preliminaries: jet bundles}

In this section notations and basic facts needed throughout
the paper are collected. The reader is referred to
\cite{AVL,KLV,VK} for further details.

Let $E$ be an $(n+m)$-dimensional manifold. The manifold of $k$-th
order jets, $k\geq0$, of $n$-dimensional submanifolds of $E$ is
denoted by $J^{k}(E,n)$ and $\pi_{k,l}:J^{k}(E,n)\longrightarrow
J^{l}(E,n)$, $k\geq l$, stands for the canonical projection. If
$E$ is fibered by a map $\pi:E\rightarrow M$ over an
$n$-dimensional manifold $M$, then $J^{k}\pi$ denotes the $k$-th
order jet manifold of local sections of $\pi$. $J^{k}\pi$ is an
open domain in $J^{k}(E,n)$. The $k-$th order jet of an
$n-$dimensional submanifold $L\subset E$ at a point $z\in L$ is
denoted by $\left[  L\right]  _{z}^{k}$. Similarly, if $\sigma$ is
a (local) section of $\pi$ and $x\in M$, then $\left[
\sigma\right]  _{x}^{k}=\left[  \sigma(U)\right]
_{\sigma(x)}^{k}$, $U$ being the domain of $\sigma$, stands for
the $k-$th order jet of $\sigma$ at $x$. The correspondence
$z\mapsto\left[  L\right]  _{z}^{k}$ defines the $k-$th
\emph{lift} of $L$
\[%
\begin{array}
[c]{cccc}%
j_{k}L: & L & \longrightarrow & J^{k}(E,n).
\end{array}
\]
Similarly, the $k-$th \emph{lift} of a local section $\sigma$ of
$\pi$
\[
j_{k}\sigma:U\rightarrow J^{k}\pi
\]
sends $x\in U$ to $[\sigma]_{x}^{k}$,  i.e.,
$j_{k}\sigma=j_{k}(\sigma (U))\circ\sigma$. Put
\[
L^{(k)}=\mathop{\rm Im}(j_{k}L),\qquad\qquad M_{\sigma}^{k}%
=\mathop{\rm Im}(j_{k}\sigma).
\]

Let $\theta_{k+1}=\left[  L \right]  _{z}^{k+1}$ be a point of
$J^{k+1}(E,n)$. Then the $R-plane$ associated with
$\theta_{k+1}$ is the subspace
\[
R_{\theta_{k+1}}=T_{\theta_{k}}(L^{(k)})
\]
of $T_{\theta_{k}}(J^{k}(E,n))$ with $\theta_{k}=\left[  L\right]
_{z}^{k}$. The correspondence $\theta_{k+1}\mapsto
R_{\theta_{k+1}}$ is biunique.
Put
\[
V_{\theta_{k+1}}=T_{\theta_{k}}(J^{k}(E,n))/R_{\theta_{k+1}}.
\]
and denote by
\[
pr_{\theta_{k+1}}:T_{\theta_{k}}(J^{k}(E,n))\longrightarrow V_{\theta_{k+1}}%
\]
the canonical projection. The vector bundle
\[
\nu_{k+1}:V_{(k+1)}\longrightarrow J^{k+1}(E,n),\qquad k\geq0,
\]
whose fiber over $\theta_{k+1}$, is $V_{\theta_{k+1}}$ is
naturally defined. By $\nu_{k,r}$ denote the pullback of $\nu_{k}$
via $\pi_{r,k},r\geq k$.

Let $\mathcal{C}(\theta_{k})\subset T_{\theta_{k}}(J^{k}\pi)$ be
the span of all $R-$planes at $\theta_{k}$. Then
$\theta_{k}\mapsto\mathcal{C}(\theta _{k})$ is the \emph{Cartan
distribution } on $J^{k}(E,n)$ denoted by $\mathcal{C}_{k}$. This
distribution can be alternatively defined as the
kernel of the $\nu_{k}$-valued \emph{Cartan form} $U_{k}$ on $J^{k}(E,n)$:%
\[
U_{k}(\xi)=pr_{\theta_{k}}(d_{\theta_{k}}\pi_{k,k-1}(\xi)) \in
V_{\theta_{k}}, \quad\xi\in T_{\theta_{k}}(J^{k}(E,n)).
\]

A diffeomorphism $\varphi:J^{k}(E,n)\rightarrow J^{k}(E,n)$ is
called \emph{contact} if it preserves the Cartan distribution.
Similarly, a vector field $Y$ on $J^{k}(E,n)$ is called
\emph{contact} if $\left[ Y,\mathcal{C}_{k}\right]
\subset\mathcal{C}_{k}$. A contact diffeomorphism $\varphi$
(respectively, a contact field $Y$) canonically lifts to a contact
diffeomorphism $\varphi^{(l)}$ (respectively, a contact field
$Y^{(l)}$) on $J^{k+l}(E,n)$).

Below the above constructions will be mainly used for
$n=2,m=1,k=1,2$. In this
case $\mathcal{C}_{1}$ is the canonical contact structure on $J^{1}%
(E,2), \mathop{\rm dim}E=3$, and the bundle $\nu_{1}$ is
$1$-dimensional. $\nu_{1}$ is canonically isomorphic to the bundle
whose fiber over $\theta\in J^{1}(E,2)$ is 
$T_{\theta}(J^{1}(E,2))/\mathcal{C}(\theta)$.

A vector field $X$ on $E$ defines a section
$s_{X}\in\Gamma(\nu_{1}), s_{X}(\theta)= pr_{\theta}(X)$. Since
$\nu_{1}$ is $1$-dimensional, the $\nu_{1}$-valued form $U_{1}$
can be presented as
\begin{equation}\label{UX}
U_{1}=U_{X}\cdot s_{X}, \quad U_{X}\in\Lambda^{1}(J^{1}(E,2)),
\end{equation}
in the domain where $s_{X}\neq0$.

Let $\mathcal{M}$ be a manifold supplied with a contact
distribution $\mathcal{C}$. An almost everywhere nonvanishing
differential form $U\in\Lambda^{1}(\mathcal{M})$ is called
\emph{contact} if it vanishes on $\mathcal{C}$. A vector field
$Y\in D(\mathcal{M})$ is contact iff
\begin{equation}
L_{X}(U)=\lambda U,\quad\lambda\in C^{\infty}(\mathcal{M}),
\label{infinitesimal_cont_inv}%
\end{equation}
for a contact form $U$. For instance, $U_{X}$ (see (\ref{UX})) is
a contact form.

If $X\in D(J^{1}(E,2)$ is a contact vector field, then $f=X\left.
\hskip1pt\vbox{\hbox{\vbox to .18 truecm{\vfill\hbox to .25 truecm
{\hfill\hfill}\vfill}\vrule}\hrule}\hskip1pt\right.
U_{1}\in\Gamma(\nu_{1})$ is called the \emph{generating function}
of $X$. $X$ is completely determined by $f$ and is denoted by
$X_{f}$ in order to underline this fact. If $U$ is a contact form
on a contact manifold $(\mathcal{M},\mathcal{C})$ and $X\in
D(\mathcal{M})$ is contact, then $f=X\left.
\hskip1pt\vbox{\hbox{\vbox to .18 truecm{\vfill\hbox to .25 truecm
{\hfill\hfill}\vfill}\vrule}\hrule}\hskip1pt\right.  U\in
C^{\infty }(\mathcal{M})$ is the \emph{generating function} of $X$
with respect to $U$.

Vector fields $X,Y\in D(\mathcal{M})$ belonging to $\mathcal{C}$
are called \emph{$\mathcal{C}$-orthogonal} if $\left[  X,Y\right]
$ also belongs to $\mathcal{C}$. Obviously, this is equivalent to
$dU(X,Y)=0$ for a contact form
$U$. Observe that $\mathcal{C}$-orthogonality is a $C^{\infty}(\mathcal{M}%
)$-linear property. A subdistribution $\mathcal{D}$ of
$\mathcal{C}$ is called
\emph{Lagrangian} if any two fields $X,Y\in\mathcal{D}$ are $\mathcal{C}%
$-orthogonal and $\mathcal{D}$ is not contained in another
distribution of bigger dimension possessing this property. If
$\dim\;\mathcal{M}=2n+1$, then $\dim\;\mathcal{D}=n$.

A local chart $\left(  x,y,u\right)  $ in $E$, where $\left(
x,y\right)  $ are interpreted as independent variables and $u$ as
the dependent one,  canonically extends to a local chart
\begin{equation}
\left(  x,y,u,u_{x}=p,u_{y}=q,u_{xx}=r,u_{xy}=s,u_{yy}=t\right)  \label{carta}%
\end{equation}
on $J^{2}(E,2)$. Functions $\left(  x,y,u,p,q\right)  $ form a
(standard) chart in $J^{1}(E,2)$. The local contact form
$U=U_{\partial_{u}}$ (see (\ref{UX})) in this chart reads
\[
U=du-pdx-qdy.
\]
Accordingly, in this chart the contact vector field corresponding
to the generating function $f$ reads
\begin{equation}
X_{f}=-f_{p}\partial_{x}-f_{q}\partial_{y}+(f-pf_{p}-qf_{q})\partial
_{u}+(f_{x}+pf_{u})\partial_{p}+(f_{y}+qf_{u})\partial_{q}. \label{contact_VF}%
\end{equation}

In the sequel we shall use $\mathcal{C}$ and $U$ to denote contact
distribution and contact form, respectively.

\section{Parabolic Monge-Amp\`ere equations}

\subsection{Monge-Amp\`ere equations}

Let $E$ be a $3$-dimensional manifold. A $k$-th order differential
equation
imposed on bidimensional submanifolds of $E$ is a hypersurface $\mathcal{E}%
\subset J^{k}(E,2)$. In a standard jet chart it looks as a
$k$-th order equation for one unknown function in two variables.
In the sequel we shall deal only with second order equations of
this kind. In a jet chart (\ref{carta}) on $J^{2}(E,2)$ such an
equation reads
\begin{equation}
F(x,y,u,p,q,r,s,t)=0. \label{SDE}%
\end{equation}

The standard subdivision of equations (\ref{SDE}) into hyperbolic,
parabolic and elliptic ones is intrinsically characterized by the
nature of singularities their multi-valued solutions (see
\cite{VinSing}) can have. Below we collected some elementary facts from solution singularity
theory we need in this paper.

Let $\theta_{2}\in J^{2}(E,2)$, $\theta_{1}=\pi_{2,1}(\theta_{2})$
and $F_{\theta_{1}}=\pi_{2,1}^{-1}(\theta_{1})$. Recall that an
$R$-plane at $\theta_{1}$ is a \emph{Lagrangian} plane in
$\mathcal{C}(\theta_{1})$, i.e., a bidimensional subspace
$R\subset\mathcal{C} (\theta_{1})$ such that $d\omega|_{R}=0$ for
a contact 1-form $\omega$ on $J^{1}(E,2)$. Denote by $LG(\theta_1)$ the
Lagrangian Grassmannian of $R$--planes at $\theta_1$. The correspondence 
$\theta_2\mapsto R_{\theta_2}$  identifies $F_{\theta_{1}}$ with the open 
subset in $LG(\theta_1)$ composed of $R$--planes that are transversal
to the fiber of the projection $\pi_{1,0}$. $F_{\theta_{1}}$  carries a natural 
affine structure, modeled over the vector space of quadratic forms on 
$R_{\theta_1}$.
An affine chart on $F_{\theta_{1}}$ is formed by restrictions of
$r,s,t$ to $F_{\theta_{1}}$. By an abuse of notation we shall use
$r,s,t$ for these restrictions as well. In order to put the theory of second 
order equations into a contact invariant form $J^2(E,2)$  must be completed  
by ``infinite" points
corresponding to $R$--planes which are non transversal to the fibers of $\pi_{1,0}$ 
(``Lagrangian projectivization").

Let 
$P\subset\mathcal{C}(\theta_{1})$ be a 1-dimensional subspace. Denote by
$l(P)$ the curve in $LG(\theta_1)$ composed of all  $R$-planes that contain P.
$l(P)$ is called the \emph{1-ray} corresponding to $P$. If $P$ is transversal to 
the fiber of $\pi_{1,0}$ passing through $\theta_1$, then
\begin{equation}\label{ray}
l(P)\cap F_{\theta_{1}} =\left\{  \theta\in
F_{\theta_{1}}|R_{\theta}\supset P\right\}
\end{equation}
is a straight line in $ F_{\theta_{1}}$ . Otherwise, this intersection is empty.

Let $R\in LG(\theta_1)$. The cone $\mathcal{V}_R\subset T_R( LG(\theta_1))$
is composed of tangents to $1$-rays passing through $R$. If $R=R_{\theta_2}$, we 
put $\mathcal{V}_{\theta_{2}}=\mathcal{V}_R\subset T_{\theta_2} (F_{\theta_{1}})$ 
by identifying $F_{\theta_{1}}$ with the corresponding domain in  $LG(\theta_1)$.
In terms of cartesian coordinates
$(\tilde{r},\tilde{s},\tilde{t})$ in
$T_{\theta_{2}}(F_{\theta_{1}})$ with respects to the basis
$\partial _{r}|_{\theta_{2}},
\partial_{s}|_{\theta_{2}},\partial_{t}|_{\theta_{2}}$ 
the equation of   $\mathcal{V}_{\theta_{2}}$ is $\tilde{r}\tilde{t}-\tilde{s}^{2}=0$.
The family $\mathcal{V}:R\mapsto\mathcal{V}_R$  will be called the 
\emph{ray distribution} on $LG(\theta_1)$. If $R=R_{\theta_2}$, $R$ (respectively, $\theta_2$) is the 
\emph{vertex} of $\mathcal{V}_R$ (resp., $\mathcal{V}_{\theta_2}$).

A coordinate description of these facts is given by the following

\begin{lemma}
Let $P$ be spanned by the vector
\[
w=\zeta_{1}\partial_{x}+\zeta_{2}\partial_{y}+\mu\partial_{u}+\eta_{1}%
\partial_{p}+\eta_{2}\partial_{q} \in T_{\theta_{1}}(J^{1}(E,2)),
\]
then the 1-ray $l(P)$ is described by equations
\begin{equation}
\left\{
\begin{array}
[c]{c}%
\zeta_{1}r+\zeta_{2}s=\eta_{1},\\
\zeta_{1}s+\zeta_{2}t=\eta_{2}.
\end{array}
\right.  \label{eq_l(P)}%
\end{equation}
In particular, $l(P)$ is tangent to the distribution  $\mathcal{V}$.
\end{lemma}

\begin{proof}
Put $p_{0}=p(\theta_{1}), q_{0}=q(\theta_{1}). $ Since $P\subset
\mathcal{C}(\theta_{1})$, we have
\begin{equation}
w\left.  \hskip1pt\vbox{\hbox{\vbox to .18 truecm{\vfill\hbox to
.25 truecm {\hfill\hfill}\vfill}\vrule}\hrule}\hskip1pt\right.
U=0\Leftrightarrow
\mu=\zeta_{1}p_{0}+\zeta_{2}q_{0}, \label{P_contatto}%
\end{equation}
and hence%
\[
w=\zeta_{1}(\partial_{x}+p_{0}\partial_{u})+\zeta_{2}(\partial_{y}%
+q_{0}\partial_{u})+\eta_{1}\partial_{p}+\eta_{2}\partial_{q}.
\]
Moreover, $R_{\theta_{2}}\supset P$ iff%
\[
w\left.  \hskip1pt\vbox{\hbox{\vbox to .18 truecm{\vfill\hbox to
.25 truecm
{\hfill\hfill}\vfill}\vrule}\hrule}\hskip1pt\right.  (dp-rdx-sdy)_{\theta_{2}%
}=w\left.  \hskip1pt\vbox{\hbox{\vbox to .18 truecm{\vfill\hbox to
.25 truecm
{\hfill\hfill}\vfill}\vrule}\hrule}\hskip1pt\right.  (dq-sdx-tdy)_{\theta_{2}%
}=0
\]
These relations are identical to (\ref{eq_l(P)}).

Obviously, the components of the tangent vector to $l(P)$ at
$\theta_{2}$  are
\begin{equation}
(\tilde{r},\tilde{s},\tilde{t})=\left(
\zeta_{2}^{2},-\zeta_{1}\zeta
_{2},\zeta_{1}^{2}\right)  \label{coord-cono}%
\end{equation}
and manifestly satisfy the equation
$\tilde{r}\tilde{t}-\tilde{s}^{2}=0$.
\end{proof}

Put
\[
\mathcal{E}_{\theta_{1}}=\mathcal{E}\cap F_{\theta_{1}}%
\]
An equation $\mathcal{E}$ is of \emph{principal type} if it
intersects transversally fibers of the projection $\pi_{2,1}$. In
such a case $\mathcal{E}_{\theta_{1}}$ is a bidimensional
submanifold of $F_{\theta_{1}},\forall\;\theta_{1}\in J^{1}(E,2)$.
Further on we assume $\mathcal{E}$ to be of principal type.

The\emph{\ symbol} \emph{of }$\mathcal{E}$ at
$\theta_{2}\in\mathcal{E}$ is
the bidimensional subspace%
\[
Smbl_{\theta_{2}}(\mathcal{E}):=T_{\theta_{2}}\mathcal{E}_{\theta_{1}}
\]
of $T_{\theta_{2}}\left(  F_{\theta_{1}}\right)  $.

A point $\theta_{2}\in\mathcal{E}$ is \emph{elliptic} (resp.,
\emph{parabolic}, or \emph{hyperbolic}) if
$\mathcal{V}_{\theta_{2}}$ intersects $Smbl_{\theta
_{2}}(\mathcal{E})$ in its vertex only (resp., along a line, or
along two
lines). So, if $\theta_{2}$ is parabolic, then $Smbl_{\theta_{2}}%
(\mathcal{E})$ is a plane tangent to the cone $\mathcal{V}_{\theta_{2}}$. 
In other words, in this case $\mathcal{E}_{\theta_{1}}$ is
tangent to the ray distribution on $F_{\theta_{1}}$.

\begin{definition}
\label{PE} An equation $\mathcal{E}$ is called \emph{elliptic
}(resp., \emph{parabolic}, or \emph{hyperbolic}) if all its points
are \emph{elliptic} (resp., \emph{parabolic}, or
\emph{hyperbolic}).
\end{definition}

\begin{lemma}
If a $1$-ray $l(P)$ is tangent to a parabolic equation
$\mathcal{E}$ at a point $\theta_{2}$, then
$l(P)\subset\mathcal{E}_{\theta_{1}}$.
\end{lemma}

\begin{proof}
The ray distribution on $F_{\theta_{1}}$ may be viewed as the
distribution of Monge's cones of a first order PDE for one unknown
function in two variables. As it is easy to see, this equation in
terms of coordinates $x_{1}=r,x_{2}=t,y=s$ on $F_{\theta_{1}}$ is
$\frac{\partial y}{\partial x_{1}}\cdot\frac{\partial y}{\partial
x_{2}}=\frac{1}{4}$. A banal computation then shows that
characteristics of this equation are exactly $1$-rays and hence
its solutions are ruled surfaces composed of $1$-rays.
\end{proof}

\begin{corollary}
\label{ruled} If $\mathcal{E}$ is a parabolic equation, then
$\mathcal{E}_{\theta_{1}}$ is a ruled surface in $F_{\theta_{1}}$
composed of $1$-rays.
\end{corollary}

For a parabolic equation $\mathcal{E}$ and a point $\theta_{1}\in
J^{1}(E,2)$ consider all 1-dimensional subspaces
$P\subset\mathcal{C}(\theta_{1})$ such that
$l(P)\subset\mathcal{E}_{\theta_{1}}$. This is a 1-parameter
family of
lines and, so, their union is a bidimensional conic surface $\mathcal{W}%
_{\theta_{1}}$ in $\mathcal{C}(\theta_{1})$. Then
\[
\theta_{1}\mapsto\mathcal{W}_{\theta_{1}},\quad\theta_{1}\in
J^{1}(E,2),
\]
is the \emph{Monge distribution} of $\mathcal{E}$. Integral curves
of this distribution are curves along which multivalued solutions
of $\mathcal{E}$ fold up. It is worth mentioning that tangent
planes to a surface $\mathcal{W}_{\theta_{1}}$ are all Lagrangian.
We omit the proof but note that Lagrangian planes are simplest
surfaces possessing this property.

Singularities of a given type of multivalued solutions of a PDE are
described by corresponding subsidiary equations. If $\mathcal{E}$
is a parabolic equation, then integral curves of the Monge distribution 
corresponding to $\mathcal{E}$ describe loci of foldings 
of its multivalued solutions.

Intrinsically, the class of Monge-Amp\`ere (MA) equations is
characterized by the property that these subsidiary equations are
as simple as possible. More precisely, this means that conic
surfaces $\mathcal{W}_{\theta_{1}}$s' must be geometrically
simplest. As we have already noticed, for parabolic equations the
simplest are Lagrangian planes. Hence parabolic MA equations (PMAs)
are conceptually distinguished as parabolic equations whose Monge
distributions are distributions of Lagrangian planes.  It will be
shown below that this definition coincides with the traditional
descriptive one.

Recall that, according to the traditional  point of
view, MA
equations are defined as equations of the form%
\begin{equation}
\label{MA_1}N(rt-s^{2})+Ar+Bs+Ct+D=0
\end{equation}
with $N,A,B,C$ and $D$ being some functions of variables
$x,y,u,p,q$. MA equations with $N=0$ are called
\emph{quasilinear}. $\Delta=B^{2}-4AC+4ND$ is the
\emph{discriminant} of (\ref{MA_1}).

\begin{proposition}
Equation (\ref{MA_1}) is elliptic (resp., parabolic, or
hyperbolic) if $\Delta<0\quad(resp., \Delta=0$, or $\Delta>0$).
\end{proposition}

\begin{proof}
As it is easy to see, the symbol of equation (\ref{MA_1}) at a
point
$\theta_{2}$ of coordinates $(r,s,t)$ is described by the equation%
\begin{equation}
\label{coord_symbl}N(t\tilde{r}+r\tilde{t}-2s\tilde{s})+A\tilde{r}+B\tilde
{s}+C\tilde{t}=0.
\end{equation}
with $(\tilde{r},\tilde{s},\tilde{t})$ subject to the relation
$\tilde {r}\tilde{t}-\tilde{s}^{2}=0$. Now the claim directly follows from
this relation and (\ref{coord_symbl}).
\end{proof}

Finally, we observe that all above definitions and constructions
are contact invariant.

\subsection{\label{parab_come_lagr}Parabolic Monge-Amp\`ere equations as Lagrangian distributions}

In this section it will be shown that the conceptual definition of
PMA equations coincides with the traditional one.
First of all, we have
\begin{proposition}
\label{equiv} Equation (\ref{MA_1}) is parabolic in the sense of
definition\,\ref{PE} if and only if $\Delta=0$. The Monge
distribution of parabolic equation (\ref{MA_1}) is a distribution
of Lagrangian planes, i.e., a Lagrangian distribution.
\end{proposition}
\begin{proof}
Let $\mathcal{E}$ stands for an equation (\ref{MA_1}). A simple direct
computation shows that $\mathcal{E}$ is tangent to the ray
distribution iff $\Delta=0$.
Coefficients of equation (\ref{MA_1}) may be thought as
functions on $J^{1}(E,2)$. Let $A_{0},\ldots,N_{0}$ be their
values at a point $\theta _{1}\in J^{1}(E,2)$. Then
\begin{equation}
N_{0}(t\tilde{r}+r\tilde{t}-2s\tilde{s})+A_{0}\tilde{r}+B_{0}\tilde{s}%
+C_{0}\tilde{t}=0 \label{simbolo}%
\end{equation}
is the equation of $\mathcal{E}_{\theta_{1}}$ in $F_{\theta_{1}}$.
If $N_{0}\neq0$ this equation describes a standard cone with the
vertex at the
point $\theta_{2}$ of coordinates $r=-C_{0}/N_{0},s=B_{0}/2N_{0}%
,t=-A_{0}/N_{0}$. So, $\mathcal{E}_{\theta_{1}}$ is the union of
$1$-rays $l(P)$ passing through $\theta_{2}$. By definition this
implies that $P\subset R_{\theta_{2}}$ and hence
$\mathcal{W}_{\theta_{1}}$ is the union of lines $P$
that belong to $R_{\theta_{2}}$. This shows that $\mathcal{W}_{\theta_{1}%
}=R_{\theta_{2}}$.

If $N_{0}=0$, then, taking into account (\ref{eq_l(P)}) and (\ref{coord-cono}), equations (\ref{MA_1}) and (\ref{simbolo}) can be rewritten in terms of $\zeta$'s and $\eta$'s as follows:
\begin{equation}
A_{0}(\zeta_{2}\eta_{2}-\zeta_{1}\eta_{1})-B_{0}\zeta_{1}\eta_{2}-D_{0}\zeta_{1}^{2}=0, \qquad
A_{0}\zeta_{2}^{2}-B_{0}\zeta_{1}\zeta
_{2}+C_{0}\zeta_{1}^{2}=0.
\label{sistema}%
\end{equation}
Assuming now that, say, $\zeta_{1}\neq 0$ (and hence $A\neq 0$) by (\ref{sistema}) we have:
\[
\zeta_{2}=\frac{B}{2A}\zeta_{1},\hspace{0.45in}\eta_{1}=-\frac{B}{2A}\eta
_{2}-\frac{D}{A}\zeta_{1}.
\]  
This shows that  in this case $\mathcal{W}_{\theta_{1}}$ is also the union of straight lines 
belonging to the Lagrangian plane 
\begin{equation} \label{N=0}
\left\langle\partial_x+\frac{B}{2A}\partial_y+(p+\frac{B}{2A}q)\partial_u-\frac{D}{A}\partial_p, 
\;\frac{B}{2A}\partial_p-\partial_q\right\rangle
\end{equation}
(see Lemma \ref{eq_l(P)}).
\end{proof}

From now on we shall denote by $\mathcal{D}_{\mathcal{E}}$ the
Monge distribution of a PMA equation $\mathcal{E}$.  Geometrical
meaning of the correspondence
$\mathcal{E}\mapsto\mathcal{D}_{\mathcal{E}}$ is easily extracted 
from the proof of the above Proposition.

\begin{corollary}
\label{vertex} If $N\neq0$, then the distribution $\mathcal{D}_{\mathcal{E}}$
associates with a point $\theta_{1}\in J^{1}(E,2)$ the $R$-plane
$R_{\theta_{2}}$ with $\theta_{2}$ being the vertex of the cone
$\mathcal{E}_{\theta_{1}}$. If $N=0$, this vertex is an "infinite point" of 
$F_{\theta_1}$, i.e., one that corresponds to an $R$--plane which is not transversal
to the fiber of $\pi_{1,0}$ passing through $\theta_1$.
\end{corollary}

A coordinate description of $\mathcal{D}_{\mathcal{E}}$ is as
follows.

\begin{proposition}
\label{coord}
If $N\neq0$, then
\begin{equation}
\mathcal{D}_{\mathcal{E}}=\left\langle\partial_{x}+p\partial_{u}-\dfrac{C}%
{N}\partial_{p}+\dfrac{B}{2N}\partial_{q},\;\partial_{y}+q\partial_{u}+\dfrac
{B}{2N}\partial_{p}-\dfrac{A}{N}\partial_{q}\right\rangle
. \label{parab_4}%
\end{equation}

If $N=0$ and $A\neq0$, then
\begin{equation}
\mathcal{D}_{\mathcal{E}}=\left\langle
A\partial_{x}+\frac{B}{2}\partial
_{y}+(Ap+\frac{B}{2}q)\partial_{u}-D\partial_{p},\;\frac{B}{2}\partial_{p}%
-A\partial_{q}\right\rangle. \label{quasilin_1}%
\end{equation}

If $N=A=0$, then $C\neq0$ and
\begin{equation}
\mathcal{D}_{\mathcal{E}}=\left\langle
\frac{B}{2}\partial_{x}+C\partial
_{y}+(\frac{B}{2}p+Cq)\partial_{u}-D\partial_{q},\;C\partial_{p}-\frac{B}{2}\partial_{q}%
\right\rangle. \label{quasilin_1}%
\end{equation}

\end{proposition}

\begin{proof}
If $(x_{0},\dots,t_{0})$ are coordinates of $\theta_{2}\in
J^{2}(E,2)$, then the $R$-plane $R_{\theta_{2}}$ is generated by
vectors
\[
\partial_{x}+p_{0}\partial_{u}+r_{0}\partial_{p}+s_{0}\partial_{q}%
\qquad\mathrm{and}\qquad\partial_{y}+q_{0}\partial_{u}+s_{0}\partial
_{p}+t_{0}\partial_{q}%
\]
at the point $\theta_{1}$. By Corollary \ref{vertex} one gets the
desired result for $N\neq0$ just by substituting coordinates of
of the vertex $\theta_{2}$ of the
cone $\mathcal{E}_{\theta_{1}}$ (see the proof of Proposition
\ref{equiv}) for $r_0, s_0, t_0$ in these expressions.

For $N=0$ the result is already given by (\ref{N=0}).
\end{proof}

In its turn, the distribution $\mathcal{D}_{\mathcal{E}}$
completely determines the equation $\mathcal{E}$. Namely, we
have

\begin{proposition}
\label{cond_PMA}Let $\mathcal{D}$ be a Lagrangian subdistribution
of $\mathcal{C}$. Then the submanifold
\[
\mathcal{E}_{\mathcal{D}}=\{\theta_{2}\in J^{2}(E,2):\dim(R_{\theta_{2}}%
\cap\mathcal{D}(\theta_{1}))>0\}
\]
of $J^{2}(E,2)$ is a parabolic Monge-Amp\`ere equation and $\mathcal{D}%
=\mathcal{D}_{\mathcal{E}_{\mathcal{D}}}$.
\end{proposition}

\begin{proof}
There is the only one point $\theta_{2}\in\left(  \mathcal{E}_{\mathcal{D}%
}\right)  _{\theta_{1}}$ such that
$R_{\theta_{2}}=\mathcal{D}(\theta_{1})$. With this exception
$\dim(R_{\theta_{2}}\cap\mathcal{D}(\theta_{1}))=1$. Hence $\left(
\mathcal{E}_{\mathcal{D}}\right)  _{\theta_{1}}$ is the union of
all $1$-rays $l(P)$ such that $P\subset\mathcal{D}(\theta_{1})$.
They all pass through the exceptional point $\theta_{2}$ and hence
constitute a cone in $F_{\theta_{1}}$. But cones composed of
$1$-rays are tangent to the ray distribution on
$F_{\theta_{1}}$ and, so, all their points are parabolic. Finally,
the last assertion directly follows from Corollary \ref{vertex}.
\end{proof}

Results of this section are summed up in the following assertion
which is the starting point of our subsequent discussion of
parabolic Monge-Amp\`ere equations.

\begin{theorem}
\label{MAC} The correspondence $D\longmapsto\mathcal{E}_{D}$
between Lagrangian distributions on $J^{1}(E,2)$ and parabolic
Monge-Amp\`ere equations is one-to-one.
\end{theorem}

The meaning of this Theorem is that it decodes the geometrical
problem hidden under analytical condition (\ref{MA_1}). Namely,
this problem is to find Legendrian submanifolds $S$ of a
$5$-dimensional contact manifold $(\mathcal{M},\mathcal{C})$ that
intersect a given Lagrangian distribution
$\mathcal{D}\subset\mathcal{C}$ in a nontrivial manner, i.e.,
\[
\dim\{T_{\theta}(S)\cap\mathcal{D}(\theta)\}>0,\,\,\forall\;\theta\in
S.
\]
The triple $\mathcal{E}=(\mathcal{M},\mathcal{C},\mathcal{D})$
encodes this problem. By this reason, in the rest of this paper
the term \textquotedblleft parabolic Monge-Amp\`ere
equation\textquotedblright\ will refer to such a triple. In
particular, equivalence and classification problems for PMAs are
interpreted as the corresponding problems for Lagrangian distributions on
5-dimensional contact manifolds.

\section{Geometry of Lagrangian distributions}

In this section we deduce some basic facts about Lagrangian
distributions on 5-dimensional contact manifolds which allow one to
reveal four natural classes of them. Throughout this section
$(\mathcal{M},\mathcal{C})$ refers to the considered contact manifold
and $\mathcal{D}$ for a Lagrangian distribution on it.

First of all, Lagrangian distributions are subdivided into
\emph{integrable} and \emph{non-integrable} ones. Accordingly, the
corresponding parabolic Monge-Amp\`ere equations are called
\emph{integrable}, or \emph{non-integrable}. It will be shown 
in the next section that all integrable PMAs are locally
contact equivalent to the equation $u_{xx}=0$. By this reason we shall
concentrate on nonintegrable PMAs.

If $\mathcal{D}$ is nonintegrable, then its \emph{first
prolongation} $\mathcal{D}_{(1)}$, i.e., the span of all vector
fields belonging to $\mathcal{D}$ and their commutators, is
3-dimensional. Moreover, we have

\begin{lemma}
\rule{0.34in}{0in}

\begin{enumerate}
\item $\mathcal{D}_{(1)}\subset\mathcal{C}$;

\item the $\mathcal{C}$-orthogonal complement $\mathcal{R}$ of $\mathcal{D}%
_{(1)}$ is a $1$-dimensional subdistribution of $\mathcal{D}$.
\end{enumerate}
\end{lemma}

\begin{proof}
(1) If (locally) $\mathcal{D}=\left\langle X,Y\right\rangle $,
then (locally) $\mathcal{D}_{(1)}=\left\langle
X,Y,[X,Y]\right\rangle $. But, by definition of
$\mathcal{C}$-orthogonality, $[X,Y]\in\mathcal{C}$.

(2) The form $dU$ restricted to $\mathcal{C}$ is nondegenerate.
So, the assertion follows from the fact that in a symplectic
linear space a Lagrangian subspace contained in a hyperplane contains, in its turn, the
skew-orthogonal complement of this hyperplane.
\end{proof}

The $1$-dimensional distribution $\mathcal{R}$ will be called
\emph{the directing distribution of }$\mathcal{D}$ (alternatively,
of $\mathcal{E}$).
This way one gets the following flag of distributions
\[
\mathcal{R}\subset\mathcal{D}\subset\mathcal{D}_{(1)}\subset\mathcal{C}.
\]

The directing distribution $\mathcal{R}$ uniquely deftermines
$\mathcal{D}_{(1)}$, which is its $\mathcal{C}$-orthogonal
complement, as well as the distribution
\[
\mathcal{D}^{\prime}=\left\{  X\in\mathcal{D}_{(1)}|\left[  X,\mathcal{R}%
\right]  \subset{\mathcal{D}}_{(1)}\right\}  .
\]

Since $\left[  X,\mathcal{R}\right]  \subset{\mathcal{D}}_{(1)}$ for any
$X\in
\mathcal{D}$, we see that $\mathcal{D}\subseteq\mathcal{D}^{\prime}%
\subseteq\mathcal{D}_{(1)}$. So, by obvious dimension arguments,
only one of the following two possibilities may  locally occur:
either $\mathcal{D}^{\prime }=\mathcal{D}$, or
$\mathcal{D}^{\prime}=\mathcal{D}_{(1)}$.  In
the first case the (non-integrable)
Lagrangian distribution $\mathcal{D}$ is called \emph{generic} and 
\emph{special} in the second one. Accordingly, the
corresponding PMAs are called generic, or special.
The following assertion directly follows from the above
definitions.

\begin{proposition}
\label{generic-MA} A generic Lagrangian distribution $\mathcal{D}$
is completely determined by its directing distribution
$\mathcal{R}$. It is no longer so for a special distribution and
in this case $\mathcal{R}$ is the characteristic distribution of
$\mathcal{D}_{(1)}$.
\end{proposition}

In this paper we shall concentrate on generic PMA equations with a
special attention to the equivalence problem. In view of
Proposition \ref{generic-MA} this problem takes part of the
equivalence problem for 1-dimensional subdistributions of
$\mathcal{C}$. By this reason we need a criterion which
distinguishes directing distributions of generic PMA equations from
other $1$-dimensional subdistributions of $\mathcal{C}$. The
notion of \emph{type} of a $1$-dimensional
subdistribution of $\mathcal{C}$, introduced below, solves this problem.

Fix a 1-dimensional distribution $\mathcal{S}\subset\mathcal{C}$
and put
\[
\mathcal{S}_{r}=\left\{  Y\in\mathcal{C}\;|\;X^{s}(Y)\in\mathcal{C}%
,\forall\;X\in\mathcal{S}\;\mathrm{and}\;0\leq s\leq r\right\}
.
\]
In the following lemma we list without a proof some obvious
properties of $\mathcal{S}_{r}$.

\begin{lemma}
\label{simple}\rule{0.34in}{0in}

\begin{itemize}
\item $\mathcal{S}_{r+1}\subset\mathcal{S}_{r}$ and $\mathcal{S}%
\subset\mathcal{S}_{r}$;

\item If (locally) $\mathcal{S}=\left\langle X\right\rangle $,
then (locally)
\[
\mathcal{S}_{r}=\left\{
Y\in\mathcal{C}\;|\;X^{s}(Y)\in\mathcal{C}, \;0\leq s\leq
r\right\}
\]

\item $\mathcal{S}_{r}$ is a finitely generated
$C^{\infty}(\mathcal{M})$-module;

\item If (locally)
$\mathcal{C}=\mathrm{Ann}(U),\;U\in\Lambda^{1}(\mathcal{M})$, and (locally)
$\mathcal{S}=\left\langle X\right\rangle $, then
\[
\mathcal{S}_{r}=\mathrm{Ann}(U,X(U),...,X^{r}(U)).
\]

\end{itemize}
\end{lemma}

Let $\mathcal{T}$ be a $C^{\infty}(\mathcal{M})$-module. An open
domain $B\subset\mathcal{M}$ is called \emph{regular} for
$\mathcal{T}$ if the localization of $\mathcal{T}$ to $B$ is a
projective $C^{\infty}(B)$-module. If $\mathcal{T}$ is finitely
generated (see \cite{nestruev}), then this localization is
isomorphic to the $C^{\infty}(B)$-module of smooth sections of a
finite dimensional vector bundle over $B$. In such a case the
dimension of this bundle is called the \emph{rank} of
$\mathcal{T}$ on $B$ and denoted by
$\mathrm{{rank}_{B}\mathcal{T}}$. Moreover, the manifold
$\mathcal{M}$ is subdivided into a number of open domains which are
regular for $\mathcal{T}$ and the set of its \emph{singular
points} which is closed and thin. In particular, vector bundles
representing localizations of $\mathcal{S}_{r}$ to its regular
domains are distributions contained in $\mathcal{C}$ and
containing $\mathcal{S}$ (Lemma \ref{simple}).

\begin{lemma}
\label{rank} Let $B\subset\mathcal{M}$ be a common regular domain
for modules
$\mathcal{S}_{r}$ and $\mathcal{S}_{r+1}$. Then either $\mathrm{{rank}%
_{B}\mathcal{S}_{r}={rank}_{B}\mathcal{S}_{r+1}+1}$, or $\mathrm{{rank}%
_{B}\mathcal{S}_{r}={rank}_{B}\mathcal{S}_{r+1}}$. In the latter
case $B$ is
regular for $\mathcal{S}_{p},p\geq r$, and $\mathrm{{rank}_{B}\mathcal{S}%
_{p}={rank}_{B}\mathcal{S}_{r}}$.
\end{lemma}

\begin{proof}
The first alternative takes place iff $X^{r+1}(U)$ is $C^{\infty}%
(\mathcal{M})$-independent of $U,X(U),\dots,X^{r}(U)$ as it is
easily seen from the last assertion of Lemma \ref{simple}. In the
second case the equality
of ranks implies that localizations of $\mathcal{S}_{r}$ and $\mathcal{S}%
_{r+1}$ to $B$ coincide. This shows that the localization of
$\mathcal{S}_{p}$ to $B$ stabilizes by starting from $p=r$.
\end{proof}

\begin{corollary}
With the exception of a thin closed set the manifold $\mathcal{M}$
is subdivided into open domains which are  regular for all
$\mathcal{S}_{p},p\geq0$. Moreover, in such a domain $B$
$\mathrm{{rank}_{B}} \mathcal{S}_{p}=4-p$, if $p\leq r$, and
$\mathrm{{rank}_{B}}\mathcal{S}_{p}=4-r$, if $p\geq r$, for an
integer $r=r(B), \,1\leq r\leq 3$.
\end{corollary}

\begin{proof}
It immediately follows from Lemma \ref{rank} that the function
$p\mapsto \mathrm{{rank}_{B}\mathcal{S}_{p}}$ monotonically decreases
up to the moment, say $p=r$, when the equality
$\mathrm{{rank}_{B}\mathcal{S}_{r}={rank}_{B}\mathcal{S}_{r+1}}$
occurs for the first time, and then stabilizes. Since $\mathcal{S}%
\subset\mathcal{S}_{p}$, the last assertion of Lemma \ref{simple}
shows that this stabilization happens at most for $p=3$, i.e., that
$r\leq3$. On the other hand, forms $U$ and $X(U)$ are independent.
Indeed, the equality $X(U)=\lambda U$ for a function $\lambda$ means that $X$ is
a contact field with vanishing generating function $f=i_{X}(U)\equiv0$. But
this implies that $X=0$ (see \ref{contact_VF}). Hence $r\geq1$.

Finally, observe that regular domains $\mathcal{S}_{p+1}$ are
obtained from those of $\mathcal{S}_{p}$ by removing from latters
some thin subsets. Since, as we have seen before, the situation
stabilizes after at most four steps by starting from $p=0$, 
existence of common regular domains for all $\mathcal{S}_{p}$'s
is guaranteed. Moreover, the above arguments shows that their union 
is everywhere dense in $\mathcal{M}$ 
\end{proof}

\begin{definition}
If $\mathcal{M}$ is the only regular domain for all
$\mathcal{S}_{p}$'s, then $\mathcal{S}$ is called \emph{regular}
on $\mathcal{M}$ and the integer $r$ by starting from which 
$\mathrm{{rank}_{B}\mathcal{S}_{r}}$ stabilizes is called the
\emph{type} of $\mathcal{S}$ and also of $X$, if
$\mathcal{S}=\left\langle X\right\rangle $.
\end{definition}

Thus Lemma \ref{rank} tells that the type of $\mathcal{S}$ can be
one of the numbers $1,2$, or $3$ only, and that the union of
domains in which $\mathcal{S}$ is regular is everywhere dense in
$\mathcal{M}$.

It is not difficult to exhibit vector fields of each these three
types. For instance, vector fields of the form $X_{f}-fX_{1}$ on
$\mathcal{M}=J^{1}(E,2)$ with everywhere non-vanishing function
$f$ are of type 1. Fields $\partial
_{x}+p\partial_{u}+q\partial_{p}$ and $\partial_{x}+p\partial_{u}%
+(xy+q)\partial_{p}$ are of type $2$ and $3$, respectively.

Now we can characterize directing distributions of generic PMA
equations

\begin{proposition}
\label{caratt} A $1$-dimensional regular distribution $\mathcal{S}%
\subset\mathcal{C}$ is the directing distribution of a generic PMA
equation if and only if it is of type $3$.
\end{proposition}

\begin{proof}
Let $\mathcal{S}$ be a regular distribution of type $3$ and
(locally) $\mathcal{S}=<X>$. Then by definition
$\mathcal{S}_{1}=\mathcal{S}^{\bot}$ and
$\mathcal{D}=\mathcal{S}_{2}$ is a bidimensional subdistribution
of $\mathcal{C}$. $\mathcal{D}$ is Lagrangian, since
$\mathcal{D}\subset \mathcal{S}_{1}=\mathcal{S}^{\bot}$ and
$\mathcal{S}\subset\mathcal{D}$. If (locally)
$\mathcal{D}=\left\langle X,Y\right\rangle $, then $\left[
X,Y\right]  \notin\mathcal{D}$. Indeed, assuming that $\left[
X,Y\right]
\in\mathcal{D}$ one finds that $X^{p}(Y)\in\mathcal{D}\subset\mathcal{C}%
,\forall p$, and hence $Y\in\mathcal{S}_{p},\forall p$. In
particular, this
implies that $\mathcal{D}=\left\langle X,Y\right\rangle \subset\mathcal{S}%
_{3}$ in contradiction with the fact that
$\mathrm{{dim}\;\mathcal{S}_{3}=1}$ if $\mathcal{S}$ is of type 3.
So, $\mathrm{{dim}\left\langle X,Y,\left[ X,Y\right]
\right\rangle =3}$, and
$\left\langle X,Y,\left[  X,Y\right]  \right\rangle \subset\mathcal{S}%
^{\bot}$, since $X(Y),X^{2}(Y)\in\mathcal{C}$. Now, by dimension
arguments, we conclude that $\left\langle X,Y,\left[  X,Y\right]
\right\rangle =\mathcal{S}^{\bot}$ and, therefore, $\mathcal{S}$
is the directing distribution of $\mathcal{D}=\mathcal{S}_{2}$.

Conversely, assume that (locally) $\mathcal{R}=\left\langle
X\right\rangle $ is the directing distribution of a generic PMA
equation corresponding to the Lagrangian distribution
$\mathcal{D}=\left\langle X,Y\right\rangle $ (locally). Then, by
definition, $\mathcal{R}_{1}=\mathcal{R}^{\bot}$.
Moreover, $\mathcal{D}=\mathcal{D}^{\prime}$ implies that $X^{2}%
(Y)\notin\mathcal{R}^{\bot}$ and hence fields $X,Y,X(Y),X^{2}(Y)$
form a local basis of $\mathcal{C}$. This shows that
$X^{3}(Y)\notin\mathcal{C}$. Indeed, the assumption
$X^{3}(Y)\in\mathcal{C}$ implies the inclusion $\left[  X,\mathcal{C}\right]
\subset\mathcal{C}$, i.e., that $X$ is a nonzero contact field
with zero generating function. Hence
$Y\notin\mathcal{R}_{3}\Leftrightarrow\left[ X,Y\right]
\notin\mathcal{R}_{2}$. From one side, this shows that
$\mathcal{D}=\mathcal{R}_{2}$ and, from other side, that $\mathcal{R}_{2}%
\neq\mathcal{R}_{3}$. So, $\mathcal{R}_{3}=\left\langle
\mathcal{R} \right\rangle
\Rightarrow\mathrm{{rank}\;\mathcal{R}=3}$.
\end{proof}

In the sequel $Z$ will denote a generator of the directing distribution of the
considered PMA $\mathcal{E}$. The following coordinate description of $\mathcal{D}=\mathcal{D}%
_{\mathcal{E}}, \,\mathcal{D}_{1}$ and $Z$ is easily obtained by a
direct computation on the basis of Proposition\,\ref{coord}:

\begin{proposition}
\label{quasi_linear}Let $\mathcal{E}$ be a quasilinear
nonintegrable PMA of the form (\ref{MA_1}) with $A\neq0$. By
normalizing its coefficients to $A=1$ one has:
$\mathcal{D}=\left\langle X_{1},X_{2}\right\rangle $ and
$\mathcal{D}_{(1)}=\left\langle X_{1},X_{2},X_{3}\right\rangle $
with
\[%
\begin{array}
[c]{l}%
X_{1}:=\partial_{x}+p\partial_{u}+\dfrac{B}{2}(\partial_{y}+q\partial
_{u})-D\partial_{p},\\
X_{2}:=\dfrac{B}{2}\partial_{p}-\partial_{q}\rule{0pt}{19pt},\\
X_{3}:=[X_{1},X_{2}]=M_{1}\left(
\partial_{y}+q\partial_{u}\right)
+M_{2}\partial_{p}\rule{0pt}{19pt},
\end{array}
\]
where
\[
M_{1}=-\dfrac{1}{2}X_{2}(B),\qquad
M_{2}=\dfrac{1}{2}X_{1}(B)+X_{2}(D),
\]
and $\mathcal{R}=\left\langle Z\right\rangle $ with
\[
Z:=M_{1}X_{1}-M_{2}X_{2}.
\]

\end{proposition}

\begin{proposition}
Let $\mathcal{E}$ be a nonintegrable PMA of the form (\ref{MA_1})
with
$N\neq0$. By normalizing its coefficients to $N=1$ one has: $\mathcal{D}%
=\left\langle X_{1},X_{2}\right\rangle $ and
$\mathcal{D}_{(1)}=\left\langle X_{1},X_{2},X_{3}\right\rangle $
with
\[%
\begin{array}
[c]{l}%
X_{1}:=\partial_{x}+p\partial_{u}-C\partial_{p}+\dfrac{B}{2}\partial_{q},\\
X_{2}:=\partial_{y}+q\partial_{u}+\dfrac{B}{2}\partial_{p}-A\partial
_{q}\rule{0pt}{19pt},\\
X_{3}:=[X_{1},X_{2}]=M_{1}\partial_{q}+M_{2}\partial_{p}\rule{0pt}{19pt},
\end{array}
\]
where
\[
M_{1}=-X_{1}(A)-\dfrac{1}{2}X_{2}(B),\qquad M_{2}=\dfrac{1}{2}X_{1}%
(B)+X_{2}(C),
\]
and $\mathcal{R}=\left\langle Z\right\rangle $ with%
\[
Z:=M_{1}X_{1}-M_{2}X_{2}.
\]

\end{proposition}

\section{\label{classif_integr}Classification of integrable parabolic Monge-Amp\`ere equtions}

For completeness we shall prove here the following, essentially
known, result in a manner that illustrate the idea of our further
approach.

\begin{theorem}
\label{INT} With the exception of singular points all integrable
parabolic Monge-Amp\`ere equations are locally contact equivalent to each other and, 
in particular, to the equation $u_{xx}=0$.
\end{theorem}

\begin{proof}
Let $\mathcal{E}=\left(
\mathcal{M},\mathcal{C},\mathcal{D}\right)  $ be an integrable PMA
equation. This means that the Lagrangian distribution $\mathcal{D}$ is
integrable and hence define a $2$-dimensional Legendrian
foliation of $\mathcal{M}$. Locally this foliation can be viewed
as a fibre bundle $\Pi:\mathcal{M}\rightarrow W$ over a
$3$-dimensional manifold $W$.  Since
$\mathrm{{ker}\,d_{\theta}(\Pi
)=\mathcal{D}(\theta)}$, the differential
$d_{\theta}(\Pi):T_{\theta}\mathcal{M} \rightarrow T_{y}W,
y=\Pi(\theta)$ sends $\mathcal{C}(\theta)$ to a bidimensional
subspace $P_{\theta}\subset T_{y}W$. This way one gets the map $\Pi_{y}:\Pi^{-1}%
(y)\rightarrow G_{3,2}(y), \,\theta\mapsto P_{\theta}$, where
$G_{3,2}(y)$ is the Grassmanian of $2$-dimensional subspaces in
$T_{y}W$. Note that $\mathrm{{dim}\,\Pi^{-1}(y)=
{dim}\,G_{3,2}(y)=2}$ and, so, the local rank of $\Pi_{y}$ may
vary from $0$ to $2$. We shall show that, with the exception of a
thin set of singular points, $\Pi_{y}$'s are of rank $2$, i.e.,
$\Pi_{y}$'s are local diffeomorphisms.

First, assume that this rank is zero for all $y\in W$, i.e.,
$\Pi_{y}$'s are locally constant maps. In this case
$\mathcal{P}_{\theta}$ does not depend on $\theta\in\Pi^{-1}(y)$
and we can put $\mathcal{P}(y)=P_{\theta}$, for a
$\theta\in\Pi^{-1}(y)$. Hence $y\mapsto\mathcal{P}(y)$ is a
distribution on $W$ and $\mathcal{C}$ is its pullback via $\Pi$.
This shows that the distribution tangent to fibers of $\Pi$, i.e.,
$\mathcal{D}$, is characteristic for $\mathcal{C}$. But a contact
distribution does not admit nonzero characteristics.

Second, if the rank of $\Pi_{y}$'s equals to one for all $y\in W$,
then
$\mathcal{M}$ is foliated by curves%
\[
\gamma_{P}=\{\theta\in\Pi^{-1}(y)|d_{\theta}\Pi(\mathcal{C}(\theta))=P\}
\]
with $P$ being a bidimensional subspace of $T_{y}W$. Locally, this
foliation may be seen as a fibre bundle
$\Pi_{0}:\mathcal{M}\rightarrow N$ over a 4-dimensional manifold
N, and $\Pi$ factorizes into the composition
\[
\mathcal{C}\overset{\Pi_{0}}{\rightarrow}N\overset{\Pi_{1}}{\rightarrow}W
\]
with $\Pi_{1}$ uniquely defined by $\Pi$ and $\Pi_{0}$. By
construction the $3$-dimensional subspace
$d_{\theta}\Pi_{0}(\mathcal{C}(\theta))\subset
T_{z}N, \,z=\Pi_{0}(\theta)$, does not depend on $\theta\in\Pi^{-1}%
_{0}(z)=\gamma_{P}$ and one can put $\mathcal{Q}(z)= d_{\theta}\Pi
_{0}(\mathcal{C}(\theta))$ for a $\theta\in\Pi^{-1}_{0}(z)$. As
before we see that $\mathcal{C}$ is the pullback via $\Pi_{0}$ of
the $3$-dimensional distribution $z\mapsto\mathcal{Q}(z)$ in
contradiction with the fact that $\mathcal{C}$ does not admit
nonzero characteristics.

Thus, except singular points, $\Pi$ is of rank $2$ and hence a
local diffeomorphism. So, locally, $\Pi_{y}$ identifies
$\Pi^{-1}(y)$ and an open domain in $G_{3,2}(y)$. By observing
that $G_{3,2}(y)=\pi_{1,0}^{-1}(y)$, with
$\pi_{1,0}:J^{1}(W,2)\rightarrow W$ being a natural projection,
one gets a local identification of $\mathcal{M}$ with an open
domain in $J^{1}(W,2)$. It is easy to see that this identification
is a contact diffeomorphism. In other words, we have proven that
any integrable Lagrangian distribution on a $5$-dimensional
contact manifold is locally equivalent to the distribution of
tangent planes to fibers of the projection $J^{1}(\mathbb{R}^{3}%
,2)\rightarrow\mathbb{R}^{3}$.

Finally, we observe that $\mathcal{D}_{\mathcal{E}}=\left\langle
\partial_{x},
\partial_{y}\right\rangle $ for the equation $\mathcal{E}=\left\{
u_{xx}=0\right\}  $ and hence this equation is integrable.
\end{proof}

\section{Projective curve bundles and  generic parabolic Monge-Amp\`ere equations}

In this section non-integrable generic Lagrangian distributions
and, therefore, the corresponding PMAs are represented as
$4$-parameter families of curves in the projective $3$-space or,
more exactly, as \emph{projective curve bundles}. Differential
invariants of single curves composing such a bundle (say,
projective curvature, torsion,etc) put together give differential
invariants of the whole bundle and consequently of the
corresponding PMA. This basic geometric idea is developed in
details in the subsequent section.

Let $N$ be a $4$-dimensional manifold. Denote by $PT_{a}^{\ast}N$
the $3$-dimensional projective space of all $1$-dimensional
subspaces of the cotangent space $T_{a}^{\ast}N$ at the
point $a\in N$. The \emph{projectivization}
$p\tau^{\ast}:PT^{\ast}N\rightarrow N$ of the cotangent
bundle $\tau^{\ast}:T^{\ast}N\rightarrow N$ is the bundle whose
total space is $PT^{\ast}N={\bigcup}_{a\in N}PT_{a}^{\ast}N$ and
the
fiber over $a\in N$ is $PT_{a}^{\ast}N$, i.e., $(p\tau^{\ast})^{-1}%
(a)=PT_{a}^{\ast}N$. A \emph{projective curve bundle} (PCB) over
$N$ is a $1$-dimensional subbundle $\pi:K\rightarrow N$ of
$p\tau^{\ast}$:
\[%
\begin{array}
[c]{ccc}%
\quad K & \hookrightarrow & PT^{\ast}N\quad\\
\pi\downarrow &  & \downarrow p\tau^{\ast}\\
\quad N & \overset{id}{\longrightarrow} & N\quad\quad
\end{array}
\]
The fiber $\pi^{-1}(y),\,y\in N$, is a smooth curve in the
projective space $PT_{y}^{\ast}N$. A diffeomorphism
$\Phi:N\rightarrow N^{\prime}$ canonically  lifts to a
diffeomorphism $PT^{\ast}N\rightarrow PT^{\ast}N^{\prime}$. This
lift sends a PCB $\pi$ over $N$ to a PCB over $N^{\prime}$.
Denote it by $\Phi\pi$. A PCB $\pi$ over $N$ and a PCB
$\pi^{\prime}$ over $N^{\prime}$ are \emph{equivalent} if there
exist a diffeomorphism $\Phi:N\rightarrow N^{\prime}$ such that
$\pi^{\prime}=\Phi\pi$.

Let $\pi:K\rightarrow N$ be a PCB and $\theta=<\rho>\in K$ with
$\rho\in T_{\pi(\theta)}^{\ast}N$. Denote by $W_{\theta}$ the
$3$-dimensional subspace
of $T_{\pi(\theta)}N$ annihilated by $\theta$, i.e.,%
\[
W_{\theta}=\left\{  \xi\in
T_{\pi(\theta)}N\,|\,\rho(\xi)=0\right\}  .
\]
Two distributions are canonically defined on $K$. First of them is
the $1$-dimensional distribution $\mathcal{R}_{\pi}$ formed by all
vertical with respect to $\pi$ vectors. The second one, denoted by
$\mathcal{C}_{\pi}$, is defined by
\[
(\mathcal{C}_{\pi})_{\theta}=\left\{  \eta\in T_{\theta}K\,|\,d_{\theta}%
\pi(\eta)\in W_{\theta}\right\}  .
\]

Obviously, $\rm{dim}\,\mathcal{C}_{\pi}=4$ and
$\mathcal{R}_{\pi}\subset\mathcal{C}_{\pi}$. If, locally,
$\mathcal{R}_{\pi}=\left\langle Z\right\rangle ,\,Z\in D(K)$, and
$\mathcal{C}_{\pi}=\mathrm{Ann}(U_{\pi}),\,U_{\pi}\in\Lambda^{1}(K)$, then
the \emph{osculating distributions} of $\pi$ are defined as
\[
\mathcal{Z}_{s}^{\pi}=\mathrm{Ann}(U_{\pi},Z(U_{\pi}),...,Z^{s}(U_{\pi})),\qquad
s=0,1,2,3.
\]
Clearly this definition does not depend on the
choice of $Z$ and $U_{\pi}$. Note that
$\mathcal{C}_{\pi}=\mathcal{Z}_{0}^{\pi}$. Also,
$\mathcal{R}_{\pi}\subset\mathcal{Z}_{s}^{\pi},\,\forall\,s\geq0$,
as it easily follows from $\left[  i_{Z},L_{Z}\right]  =0$ and
$U_{\pi}(Z)=0$. Moreover,  forms
$U_{\pi},Z(U_{\pi}),...,Z^{3}(U_{\pi})$ are, generically, independent and, so, by
dimension arguments, $\mathcal{R}_{\pi}=\mathcal{Z}_{3}^{\pi}$.

We say that $\pi$ is a \emph{regular PCB} iff the following two
conditions are satisfied: (i)
$\mathcal{R}_{\pi}=\mathcal{Z}_{3}^{\pi}$ and (ii)
$\mathcal{C}_{\pi}$ is a contact structure on $K$. We emphasize
that regularity is a generic condition. Moreover, conditions
(i)-(ii) are equivalent to the fact that $\mathcal{R}_{\pi}$ is of
type 3 with respect to the contact distribution
$\mathcal{C}_{\pi}$. So, by Proposition \ref{caratt}, the
distribution
\[
\mathcal{D}_{\pi}=\{X\in\mathcal{R}_{\pi}^{\mathcal{\bot}}|L_{X}%
(\mathcal{R}_{\pi})\subset\mathcal{R}_{\pi}^{\mathcal{\bot}}\}.
\]
with $\mathcal{R}_{\pi}^{\bot}$ being the
$\mathcal{C}_{\pi}$-orthogonal complement of $\mathcal{R}_{\pi}$
is bidimensional and Lagrangian for a regular PCB $\pi$. Thus we
have

\begin{theorem}
\label{PCB_equiv_PMA}\label{GNR}If $\pi$ is a regular projective curve bundle,
 then $\mathcal{D}_{\pi}$ is a Lagrangian subdistribution of $\mathcal{C}_{\pi}$ and
$(K,\mathcal{C}_{\pi},\mathcal{D}_{\pi})$ is a generic parabolic 
Monge-Amp\`ere equation whose
directing distribution is $\mathcal{R}_{\pi}$. Conversely, a
generic parabolic Monge-Amp\`ere equation locally determines a regular 
projective curve bundle.
\end{theorem}

\begin{proof}
The first assertion of the Theorem is already proved. It remains
to represent a generic PMA $\left(
\mathcal{M},\mathcal{C},\mathcal{D}\right)  $ as a regular PCB.
Integral curves of its directing distribution $\mathcal{R}$
foliate $\mathcal{M}$. Locally, this foliation may be considered
as a fiber bundle $\pi:\mathcal{M}\rightarrow N$ over a
4-dimensional manifold $N$. Since
$\mathcal{R}(\theta)\subset\mathcal{C}(\theta)$, the subspace
$V_{\theta }=d_{\theta}\pi(\mathcal{C}(\theta))\subset
T_{y}N,\,y=\pi(\theta)$, is $3$-dimensional. Put
$\Gamma_{y}=\pi^{-1}(y)$. The map $\pi_{y}:\Gamma _{y}\rightarrow
G_{4,3}(y),\,\theta\mapsto V_{\theta}$, with $G_{4,3}(y)$ being
the Grassmanian of $3$-dimensional subspaces in $T_{y}N$, is
almost everywhere of rank $1$. Indeed, the assumption that (locally)
this rank is zero leads, as in the proof of Theorem \ref{INT}, to
conclude that (locally) the contact distribution $\mathcal{C}$ is
the pullback via $\pi$ of a $3$-dimensional distribution on $N$.

The correspondence $\iota_{y}:G_{4,3}(y)\rightarrow
PT_{y}^{\ast}N$ that sends a $3$-dimensional subspace $V\subset
T_{y}N$ to $\mathrm{Ann}(V)$ is, obviously, a diffeomorphism. Hence the
composition $\iota_{y}\circ\pi_{y}$ is a local embedding 
except for a thin subset of singular points. Now, it is easy to
see that images of $\Gamma_{y}$'s via $\iota_{y}\circ\pi_{y}$'s
give the required PCB.
\end{proof}

The above construction associating a PCB with a given PMA is
manifestly
functorial, i.e., an equivalence $F:\left(  \mathcal{M},\mathcal{C}%
,\mathcal{D}\right)  \longrightarrow\left(  \mathcal{M}^{\prime}%
,\mathcal{C}^{\prime},\mathcal{D}^{\prime}\right)  $ of PMAs
induces an equivalence $\Phi:\left(  N,K,\pi\right)
\longrightarrow\left(  N^{\prime },K^{\prime},\pi^{\prime}\right)
$ of associated PCBs. Indeed, $F$ sends $\mathcal{R}$ to
$\mathcal{R}^{\prime}$ and hence integral curves of
$\mathcal{R}$ (locally, fibers of $\pi$) to integral curves of $\mathcal{R}%
^{\prime}$ (locally, fibers of $\pi^{\prime}$). This defines a map
$\Phi$ of the variety $N$ of fibers of $\pi$ to the variety
$N^{\prime}$ of fibers of $\pi^{\prime}$, etc. Thus we have

\begin{corollary}
The problem of local contact classification of generic parabolic 
Monge-Amp\`ere equations is
equivalent to the problem of local classification of regular PCBs
with respect to diffeomorphisms of base manifolds.
\end{corollary}

Now we observe that there is another, in a sense, dual PCB
associated with a given PMA equation. Namely, associate with a
point $\theta\in K$ the line
$L_{\theta}=d_{\theta}\pi(D_{\pi}(\theta))\subset
T_{y}N,\,y=\pi(\theta)$. The correspondence $\theta\mapsto
L_{\theta}\in PT_{y}N$, where $PT_{y}N$ denotes the projective
space of lines in $T_{y}N$, defines a map of $\pi^{-1}(y)$ to
$PT_{y}N$, i.e., a curve (with singularities)  in $PT_{y}N$. 
As we have seen previously, this
(locally) defines  a $1$-dimensional subbundle in the
projectivization $PTN$ of $TN$. It will be called the
\emph{second} PCB associated with the considered PMA.

PCBs may be considered as canonical models of PMAs. Besides, 
they suggest a geometrically transparent construction of scalar
differential invariants of PMAs.

Let $\mathcal{I}$ be a scalar projective differential invariant of
curves in $\mathbb{R}P^{3}$, say, the \emph{projective curvature}
(see \cite{Wilc,SuB}), $\theta\in K$ and $y=\pi(\theta)$. The
value of this invariant for the curve $\Gamma_{y}=\pi^{-1}(y)$ in
$PT_{y}^{\ast}$ is a function on $\Gamma_{y}$. Denote it by
$\mathcal{I}_{\pi,y}$ and put $\mathcal{I}_{\pi}(\theta
)=\mathcal{I}_{\pi,y}(\theta)$ if $\theta\in\Gamma_{y}\subset K$.
Then, obviously, $\mathcal{I}_{\pi}\in C^{\infty}(K)$ is a
differential invariant of the PCB $\pi$ and, therefore, of the PMA 
associated with $\pi$.

\begin{theorem}
The differential invariants of the form $\mathcal{I}_{\Psi}$ are
sufficient for a complete classification of generic parabolic Monge-Amp\`ere equations
on the basis of the "principle of $n$-invariants".
\end{theorem}

\begin{proof}
According to the "principle of $n$-invariants", it is sufficient
to construct $n=\mathrm{{dim}\,\mathcal{M}=5}$ independent
differential invariants of PMAs in order to solve the
classification problem. Such invariants of the required form will
be constructed in the next section.
\end{proof}

For further details concerning the "principle of $n$-invariants", the reader is referred to \cite{AVL,SDI}.

\section{Differential invariants of generic projective curve bundles}

Let $\pi:K\rightarrow N$ be a regular PCB and
$\mathcal{R}_{\pi}=<Z>$, $\mathcal{D}_{\pi}=<Z,X>$ and
$\mathcal{C}_{\pi }=\mathrm{Ann}(U)$. Here  $Z$ and  $U$ are unique 
up to a nowhere vanishing functional  factor,
while $X$ is unique up to a transformation $X\longmapsto
gX+\varphi Z,\,g,\varphi\in C^{\infty}(K)$ with nowhere vanishing
$g$.

Since the considered PCB is regular, we have the following flag of
distributions
\[
\mathcal{R}_{\pi}\subset\mathcal{D}_{\pi}\subset\mathcal{R}_{\pi}^{\bot
}\subset\mathcal{C}_{\pi}\subset\mathcal{D}(K)
\]
of dimensions increasing from $1$ to $5$, respectively. In terms
of $X,Z$ and $U$ they are described as follows:

\begin{proposition}
\label{theorem_1}Locally, with the exception of a thin set of
singular points we have:

\begin{enumerate}
\item $\mathcal{R}_{\pi}=\mathrm{Ann}(Z^{i}(U):i=0,1,2,3)=\left\langle
Z\right\rangle ; $

\item $\mathcal{D}_{\pi}=\mathrm{Ann}(U,Z(U),Z^{2}(U))=\left\langle
Z,X\right\rangle ;$

\item $\mathcal{R}_{\pi}^{\bot}=\mathrm{Ann}(U,Z(U))=<Z,X,Z(X)>;$

\item $\mathcal{C}_{\pi}=\mathrm{Ann}(U)=<Z,Z^{i}(X):i=0,1,2>;$

\item $\{Z,X,Z(X),Z^{2}(X),Z^{3}(X)\}$ is a base of the
$C^{\infty}(K)$-module $\mathcal{D}(K)$;

\item for some functions $r_{i}\in C^{\infty}(K)$
\begin{equation}
Z^{4}(U)+r_{1}Z^{3}(U)+r_{2}Z^{2}(U)+r_{3}Z(U)+r_{4}U=0
\label{curva_proj_esterna}%
\end{equation}

\end{enumerate}
\end{proposition}

\begin{proof}
Assertions (1)-(3) are direct consequences of Lemma \ref{simple},
Proposition \ref{caratt} and definitions. Assertion (4) follows
from the fact that $Z^{2}(X)$ is independent of  fields $Z,X,Z(X)$. Indeed,
otherwise $Z$ would be a characteristic of the distribution
$\mathcal{R}_{\pi}^{\bot}=\quad<Z,X,Z(X)>$ in contradiction with
the fact that $\mathcal{R}_{\pi}$ is of rank 3. Similarly,
$Z^{3}(X)$ is independent of $Z,X,Z(X),Z^{2}(X)$, since,
otherwise, $Z$ would be a characteristic of $\mathcal{C}_{\pi}$.
This proves (5).

Finally, forms $Z^{s}(U),\,s\geq0$, are annihilated by $Z$. Since
$\mathrm{{dim}\,\mathcal{M}=5}$ this implies that $Z^{4}(U)$
depends on $Z^{s}(U),\,s\leq 3$, i.e.,
(\ref{curva_proj_esterna}).
\end{proof}

\begin{corollary}
\label{pippo} If $0\leq k,l\leq3$ and $k+l=3$, then
\[
Z^{k}(X)\lrcorner Z^{l}(U)=(-1)^{k+1}Z^{3}(X)\lrcorner U\neq0.
\]

\end{corollary}

\begin{proof}
Assertions (5) and (6) of the above Proposition show that
$Z^{3}(X)$ completes a basis of $\mathcal{C}$ to a basis of
$D(K)$. So, $Z^{3}(X)\lrcorner U$ is a nowhere vanishing function.
Next, by applying the standard formula $\left[  i_{X},L_{Y}\right]
=i_{\left[ X,Y\right]  }$ we find
\begin{eqnarray}
Z^{k}(X)\lrcorner Z^{l}(U)  &  =\left[  i_{Z^{k}(X)},L_{Z}\right]
(Z^{l-1}(U))+L_{Z}(Z^{k}(X)\lrcorner Z^{l-1}(U))=\nonumber\label{r+s}\\
&  -Z^{k+1}(X)\lrcorner Z^{l-1}(U)+L_{Z}(Z^{k}(X)\lrcorner
Z^{l-1}(U))
\end{eqnarray}
But, according to Proposition \ref{theorem_1}, (1)-(3),
$Z^{r}(X)\lrcorner Z^{s}(U)=0$, if $r+s=2$. So, for $k+l=3$
relation (\ref{r+s}) becomes
\[
Z^{k}(X)\lrcorner Z^{l}(U)=-Z^{k+1}(X)\lrcorner Z^{l-1}(U)
\]

\end{proof}

Decompose the vector field $Z^{4}(X)$ with respect to the base  (5) of
Proposition \ref{theorem_1} 
\begin{equation}
Z^{4}(X)+\rho_{1}Z^{3}(X)+\rho_{2}Z^{2}(X)+\rho_{3}Z(X)+\rho_{4}X+\rho_{5}Z=0
\label{curva_proj_interna_*}, \quad \rho_{i}\in C^{\infty}(K).%
\end{equation}
Then we have
\begin{equation}
Z^{4}(X\wedge Z)+\rho_{1}Z^{3}(X\wedge Z)+\rho_{2}Z^{2}(X\wedge
Z)+\rho
_{3}Z(X\wedge Z)+\rho_{4}X\wedge Z=0. \label{curva_proj_interna}%
\end{equation}

The last is a constraint imposed on the bivector $X\wedge Z$.
This bivector generates the distribution $\mathcal{D}_{\pi}$ and hence  is
unique up to a functional factor.

\begin{proposition}\label{ohoho}
Functions $r_{i}$'s in decomposition (\ref{curva_proj_esterna})
are expressed in terms of iterated Lie derivatives
$Z^{i}(U),Z^{j}(X)$ as follows:
\[%
\begin{array}
[c]{l}%
r_{1}=-\dfrac{X\lrcorner Z^{4}(U)}{X\lrcorner Z^{3}(U)},\qquad r_{2}%
=-\dfrac{Z(X)\lrcorner Z^{4}(U)+r_{1}Z(X)\lrcorner
Z^{3}(U)}{Z(X)\lrcorner
Z^{2}(U)},\vspace{0.1in}\\
r_{3}=-\dfrac{Z^{2}(X)\lrcorner Z^{4}(U)+r_{1}Z^{2}(X)\lrcorner Z^{3}%
(U)+r_{2}Z^{2}(X)\lrcorner Z^{2}(U)}{Z^{2}(X)\lrcorner Z(U)},\vspace{0.1in}\\
r_{4}=-\dfrac{Z^{3}(X)\lrcorner Z^{4}(U)+r_{1}Z^{3}(X)\lrcorner Z^{3}%
(U)+r_{2}Z^{3}(X)\lrcorner Z^{2}(U)+r_{3}Z^{3}(X)\lrcorner Z(U)}%
{Z^{3}(X)\lrcorner U}.
\end{array}
\]

\end{proposition}

\begin{proof}
By subsequently inserting fields $Z^{s}(X),\,0\leq s\leq3$, into the left hand side
of (\ref{curva_proj_esterna}) one easily gets the result by taking
into account Proposition \ref{theorem_1} and Corollary
\ref{pippo}.
\end{proof}

Similarly, we have

\begin{proposition}
Functions $\rho_{i}$ in decomposition (\ref{curva_proj_interna})
are expressed in terms of iterated Lie derivatives
$Z^{i}(U),Z^{j}(X)$ as follows:
\[%
\begin{array}
[c]{l}%
\rho_{1}=-\dfrac{Z^{4}(X)\lrcorner U}{Z^{3}(X)\lrcorner
U},\qquad\rho
_{2}=-\dfrac{Z^{4}(X)\lrcorner Z(U)+\rho_{1}Z^{3}(X)\lrcorner Z(U)}%
{Z^{2}(X)\lrcorner Z(U)},\vspace{0.1in}\\
\rho_{3}=-\dfrac{Z^{4}(X)\lrcorner
Z^{2}(U)+\rho_{1}Z^{3}(X)\lrcorner
Z^{2}(U)+\rho_{2}Z^{2}(X)\lrcorner Z^{2}(U)}{Z(X)\lrcorner Z^{2}(U)}%
,\vspace{0.1in}\\
\rho_{4}=-\dfrac{Z^{4}(X)\lrcorner
Z^{3}(U)+\rho_{1}Z^{3}(X)\lrcorner
Z^{3}(U)+\rho_{2}Z^{2}(X)\lrcorner Z^{3}(U)+\rho_{3}Z(X)\lrcorner Z^{3}%
(U)}{X\lrcorner Z^{3}(U)}.
\end{array}
\]

\end{proposition}

\begin{proof}
As in the proof of Proposition \ref{ohoho} we get the desired result by subsequently inserting the
vector field in the left hand side of (\ref{curva_proj_interna_*})
into $1$-forms $Z^{s}(U), \,0\leq s\leq3$.
\end{proof}

\begin{remark}
By introducing functions $\alpha_{kl}=Z^{k}(X)\lrcorner Z^{l}(U)$
we see that $r_{i}$'s and $\rho_{i}$'s are rational functions of
$\alpha_{kl}$'s:
\[
r_{1}=-\frac{\alpha_{04}}{\alpha_{03}},\qquad r_{2}=\frac{\alpha_{04}\alpha_{13}-\alpha_{03}%
\alpha_{14}}{\alpha_{03}\alpha_{12}},\qquad
etc.
\]

\end{remark}

Differential invariants we are going to construct are projective
differential invariants of curves composing the considered PCB.
Clearly, it is not possible to describe explicitly these curves.
So, the problem is how to express these invariants in terms of the
data at our disposal, i.e., $X,Z$ and $U$. In what follows this
problem is solved on the basis of a geometrically clear analogy.
For instance, the field $Z$ restricted to one of these curves may
be thought as the derivation with respect to a parameter along
this curve. So, with this interpretation in mind it
is sufficient to mimic a known construction of projective
differential invariants for curves in order to obtain the desired
result. By following the classical book  \cite{Wilc} by Wilczynsky 
we shall use Wilczynsky's
$p_{i}$'s and $q_{j}$'s instead of  $r_{i}$'s and
$\rho_{j}$'s:
\begin{equation}
r_{1}=4p_{1},\qquad r_{2}=6p_{2},\qquad r_{3}=4p_{3},\qquad
r_{4}=p_{4}
\label{Wilc_coeff_esterna}%
\end{equation}
and%
\begin{equation}
\rho_{1}=4q_{1},\qquad\rho_{2}=6q_{2},\qquad\rho_{3}=4q_{3},\qquad\rho
_{4}=q_{4}. \label{Wilc_coeff_interna}%
\end{equation}
In these terms relations (\ref{curva_proj_esterna})
and
(\ref{curva_proj_interna}) read%
\begin{equation}
Z^{4}(U)+4p_{1}Z^{3}(U)+6p_{2}Z^{2}(U)+4p_{3}Z(U)+p_{4}U=0,
\label{curva_proj_esterna_adatt}%
\end{equation}%
\begin{equation}
Z^{4}(X\wedge Z)+4q_{1}Z^{3}(X\wedge Z)+6q_{2}Z^{2}(X\wedge
Z)+4q_{3}Z(X\wedge
Z)+q_{4}X\wedge Z=0. \label{curva_proj_interna_adatt}%
\end{equation}
Formally, they look identical to Wilczynsky's formulas (see
equation (1), page 238 of \cite{Wilc}).

It should be stressed $Z$ and $U$ are unique up to a "gauge" transformation
$(Z,U)\longmapsto
(\overline{Z},\overline{U})$%
\begin{equation}
Z=f\overline{Z},\hspace{0.3in}U=h\overline{U} \label{gauge}%
\end{equation}
with nowhere vanishing $f,h\in C^{\infty}(K)$. The corresponding
transformation of coefficients
$\{p_{i}\}\longmapsto\{\bar{p}_{i}\}$ can be
easily obtained from (\ref{curva_proj_esterna_adatt}) by a direct computation:%

\begin{equation}%
\begin{array}
[c]{l}%
p_{1}\mapsto\bar{p}_{1}={\dfrac{Z(h)}{h}}+{\dfrac{p_{1}}{f}}+\,{\dfrac
{3Z(f)}{2f},\rule[-0.2in]{0in}{0.3in}}\\
p_{2}\mapsto\bar{p}_{2}={\dfrac{p_{2}}{{f}^{2}}+\dfrac{2p_{1}\,Z(f)}{{f}%
^{2}}}+{\dfrac{2p_{1}\,Z(h)}{hf}}+3\,{\dfrac{Z(h)Z(f)}{hf}}+{\dfrac{7Z(f)^{2}%
}{6{f}^{2}}\rule[-0.16in]{0in}{0.3in}}\\
\qquad\qquad+{\dfrac{2Z^{2}(f)}{3f}}+{\dfrac{Z^{2}(h)}{h},\rule[-0.2in]%
{0in}{0.3in}}\\
p_{3}\mapsto\bar{p}_{3}={\dfrac{p_{3}}{{f}^{3}}+{\dfrac{3p_{2}Z(h)}{h{f}%
^{2}}}+{\dfrac{3p_{2}Z(f)}{2{f}^{3}}}+\dfrac{p_{1}\,Z(f)^{2}}{{f}^{3}}%
+\dfrac{p_{1}Z^{2}(f)}{{f}^{2}}\rule[-0.16in]{0in}{0.3in}}\\
\qquad\qquad{+{\dfrac{6p_{1}Z(h)Z(f)}{h{f}^{2}}}}+{\dfrac{3p_{1}Z^{2}(h)}%
{hf}+{\dfrac{Z^{3}(f)}{4f}}+{\dfrac{Z(f)^{3}}{4{f}^{3}}}+{\dfrac{Z^{3}(h)}{h}%
}\rule[-0.16in]{0in}{0.3in}}\\
\qquad\qquad{+{\dfrac{9Z(f)Z^{2}(h)}{2hf}}}+{\dfrac{7Z(h)Z(f)^{2}}{2h{f}^{2}%
}+{\dfrac{2Z(h)Z^{2}(f)}{hf}}\rule[-0.16in]{0in}{0.3in}}\\
\qquad\qquad{+{\dfrac{Z^{2}(f)Z(f)}{{f}^{2}}},\rule[-0.2in]{0in}{0.3in}}\\
p_{4}\mapsto\bar{p}_{4}={\dfrac{p_{4}}{{f}^{4}}+{\dfrac{4p_{3}Z(h)}{h{f}%
^{3}}}}+{\dfrac{6p_{2}Z^{2}(h)}{h{f}^{2}}}+{\dfrac{6p_{2}Z(h)Z(f)}{h{f}^{3}%
}\rule[-0.16in]{0in}{0.3in}}\\
\qquad\qquad{+{\dfrac{4p_{1}\,Z^{3}(h)}{hf}+\dfrac{4p_{1}Z(h)Z^{2}(f)}{h{f}^{2}%
}}}+{{\dfrac{4p_{1}Z(h)Z(f)^{2}}{h{f}^{3}}}+{\dfrac{Z^{4}(h)}{h}}\rule[-0.16in]%
{0in}{0.3in}}\\
\qquad\qquad+{\dfrac{12p_{1}Z(f)Z^{2}(h)}{h{f}^{2}}+{\dfrac{6Z(f)Z^{3}(h)}{hf}%
}+{\dfrac{7Z(f)^{2}Z^{2}(h)}{h{f}^{2}}}\rule[-0.16in]{0in}{0.3in}}\\
\qquad\qquad{+{\dfrac{Z(h)Z^{3}(f)}{hf}}+{\dfrac{4Z(h)Z^{2}(f)Z(f)}{h{f}^{2}}%
}+{\dfrac{Z(f)^{3}Z(h)}{h{f}^{3}}}\rule[-0.16in]{0in}{0.3in}}\\
\qquad\qquad+{\dfrac{4Z^{2}(h)Z^{2}(f)}{hf}.}%
\end{array}
\label{trasf_Pi}%
\end{equation}

Now the problem is to combine $p_{i}$'s in a way to obtain
expressions which are invariant with respect to transformations
(\ref{trasf_Pi}). To this end we first normalize $(Z,U)$ by the
condition $p_{1}=0$. This can be easily done
with $f=1$ and a solution $h$ of the equation%
\[
{Z(h)}+p_{1}h{=0.}%
\]
After this normalization, equation
(\ref{curva_proj_esterna_adatt}) takes a
simpler form%
\begin{equation}
Z^{4}(\overline{U})+6P_{2}Z^{2}(\overline{U})+4P_{3}Z(\overline{U}%
)+P_{4}\overline{U}=0 \label{curva_esterna_semi_can}%
\end{equation}
with
\begin{equation}%
\begin{array}
[c]{l}%
P_{2}:=p_{2}-p_{1}^{2}-Z(p_{1}),\vspace{0.13in}\\
P_{3}:=p_{3}-Z(Z(p_{1}))-3p_{1}p_{2}+2p_{1}^{3},\vspace{0.13in}\\
P_{4}:=p_{4}-4p_{1}p_{3}-3p_{1}^{4}-Z(Z(Z(p_{1})))+3Z(p_{1})^{2}\vspace
{0.1in}\\
\qquad+6p_{1}^{2}Z(p_{1})+6p_{1}^{2}p_{2}-6Z(p_{1})p_{2}.
\end{array}
\label{Ps_per_tipo_3}%
\end{equation}

\begin{proposition}
Transformations (\ref{gauge}) preserving the normalization
$p_{1}=0$ are
subject to the condition%
\begin{equation}
Z(h)+{\dfrac{3}{2}Z(\ln(f))h=0.} \label{gauge_cond}%
\end{equation}

\end{proposition}

\begin{proof}
A direct computation.
\end{proof}

Now the problem reduces to finding invariant combinations of
$P_{2},P_{3}$ and $P_{4}$ with respect to normalized, i.e.,
respecting condition (\ref{gauge_cond}), transformations
(\ref{gauge}). This can be done, for instance, by mimicking the
construction of projective curvature and torsion in \cite{Wilc}.
Namely, introduce first the functions%
\begin{equation}%
\begin{array}
[c]{l}%
\Theta_{3}=P_{3}-\frac{3}{2}Z(P_{2}),\vspace{0.1in}\\
\Theta_{4}=P_{4}-2Z(P_{3})+\frac{6}{5}Z(Z(P_{2}))-\frac{81}{25}P_{2}%
^{2},\vspace{0.1in}\\
\Theta_{3\cdot1}=6\Theta_{3}Z(Z(\Theta_{3}))-7Z(\Theta_{3})^{2}-\frac{108}%
{5}P_{2}\Theta_{3}^{2}%
\end{array}
\label{P_tipo_3}%
\end{equation}
They are \emph{semi-invariant} with respect to normalized
transformations (\ref{gauge}), i.e., they are transformed
according to formulas
\begin{equation}
\overline{\Theta}_{3}=\frac{\Theta_{3}}{f^{3}},\hspace{0.45in}\overline
{\Theta}_{4}=\frac{\Theta_{4}}{f^{4}},\hspace{0.45in}\overline{\Theta}%
_{3\cdot1}=\frac{\Theta_{3\cdot1}}{f^{8}}. \label{theta_transformation}%
\end{equation}

Obviously, the following combinations of the $\Theta_{i}$'s
\begin{equation}
\kappa_{1}=\frac{\Theta_{4}}{\Theta_{3}^{4/3}},\hspace{0.45in}\kappa_{2}%
=\frac{\Theta_{3\cdot1}}{\Theta_{3}^{8/3}} \label{first_proj_inv}%
\end{equation}
are invariant with respect to normalized transformations
(\ref{gauge}).

Thus we have

\begin{proposition}
$\kappa_{1}$ and $\kappa_{2}$ are scalar differential invariants
of parabolic Monge-Ampere equations (\ref{MA_1}) with respect to
contact transformations.
\end{proposition}

Explicit coordinate expressions of $\kappa_{1}$ and $\kappa_{2}$ in terms of
coefficients of PMA (\ref{MA_1}) can be straightforwardly computed
by using those of $Z$ and $U$. However, they are too cumbersome 
to be reported here.

Another invariant, which can be readily extracted from
(\ref{theta_transformation}), is the \emph{invariant vector field}
\begin{equation}
\label{N1}N_{1}=\Theta_{3}^{-1/3}Z.
\end{equation}

Another set of scalar differential invariants can be constructed
in a similar manner by starting from equation
(\ref{curva_proj_interna}). Indeed, the vector field $Z$ and the
bivector field $X\wedge Z$ generating distributions $\mathcal{R}$
and $\mathcal{D}$, respectively, are unique up to ``gauge" transformations
\begin{equation}
Z=f\overline{Z},\hspace{0.45in}X\wedge
Z=g\overline{X}\wedge\overline{Z}
\label{gauge_interno}%
\end{equation}
with nowhere vanishing $f,g\in C^{\infty}(K)$.

It is easy to check that coefficients $q_{i}$'s are transformed
according to formulas (\ref{trasf_Pi}) and one can repeat what was
already done previously in the case of equation
(\ref{curva_proj_esterna}). In particular, equation
(\ref{curva_proj_interna}) can be normalized in the form
\begin{equation}
Z^{4}(\overline{X}\wedge\overline{Z})+6Q_{2}Z^{2}(\overline{X}\wedge
\overline{Z})+4Q_{3}Z(\overline{X}\wedge\overline{Z})+Q_{4}\overline{X}%
\wedge\overline{Z}=0 \label{curva_proj_interna_ad}%
\end{equation}
with
\begin{equation}%
\begin{array}
[c]{l}%
Q_{2}=q_{2}-q_{1}^{2}-Z(q_{1}),\vspace{0.13in}\\
Q_{3}=q_{3}-Z(Z(q_{1}))-3q_{1}q_{2}+2q_{1}^{3},\vspace{0.13in}\\
Q_{4}=q_{4}-4q_{1}q_{3}-3q_{1}^{4}-Z(Z(Z(q_{1})))+3Z(q_{1})^{2}\vspace
{0.1in}\\
\qquad+6q_{1}^{2}Z(q_{1})+6q_{1}^{2}q_{2}-6Z(q_{1})q_{2}.
\end{array}
\label{Qs_curva_proj_int}%
\end{equation}

This way we obtain the following scalar differential invariants of
PMAs
\begin{equation}
\tau_{1}:=\frac{\Lambda_{4}}{\Lambda_{3}^{4/3}},\hspace{0.45in}\tau_{2}%
:=\frac{\Lambda_{3\cdot1}}{\Lambda_{3}^{8/3}} \label{second_inv}%
\end{equation}
where \emph{semi-invariants} $\Lambda_{i}$ are
\begin{equation}%
\begin{array}
[c]{l}%
\Lambda_{3}=Q_{3}-\frac{3}{2}Z(Q_{2}),\vspace{0.1in}\\
\Lambda_{4}=Q_{4}-2Z(Q_{3})+\frac{6}{5}Z(Z(Q_{2}))-\frac{81}{25}Q_{2}%
^{2},\vspace{0.1in}\\
\Lambda_{3\cdot1}=6\Lambda_{3}Z(Z(\Lambda_{3}))-7Z(\Lambda_{3})^{2}-\frac
{108}{5}Q_{2}\Lambda_{3}^{2}.
\end{array}
\label{Lambdas}%
\end{equation}

Similarly,
\begin{equation}
N_{2}=\Lambda_{3}^{-1/3}Z.
\end{equation}
is an invariant vector field.

\begin{remark}
The observed parallelism in construction of two sets of
differential invariants is explained by the fact that in both
cases we compute the same invariants in two different PCB, namely,
the first and the second ones, associated with the considered PMA.
\end{remark}

Since vector fields $N_{1}$ and $N_{2}$ are invariant and
$N_{1}=\lambda N_{2}$ the factor
$\lambda=\Theta_{3}^{-1/3}\Lambda_{3}^{1/3}$ is a scalar
differential invariant. So,
\[
\gamma_{3}:=\lambda^{3}=\frac{\Lambda_{3}}{\Theta_{3}}%
\]
is a scalar differential invariant of a new ``mixed'' type as well
as 
\[
\gamma_{4}:=\frac{\Lambda_{4}}{\Theta_{4}},
\qquad\lambda_{3\cdot1}:=
\frac{\Lambda_{3\cdot1}}{\Theta_{3\cdot1}}%
\]
\begin{remark}\label{teta3}
The above invariants were constructed assuming that the function $\Theta_3$
is everywhere different from zero. So, the class of PMAs with $\Theta_3\equiv 0$ requires a 
special study. Partially it consists of non-generic PMAs (see \cite{V}), partially of degenerate
generic ones. The corresponding PCBs are composed of projective curves for which curvature
and torsion are not defined.
\end{remark}

By applying invariant vector fields $N_{1}$ and $N_{2}$  to already constructed scalar differential invariants
\[
\kappa_{1}, \quad\kappa_{2}, \quad\tau_{1}, \quad\tau_{2},
\quad\gamma_{3},
\quad\gamma_{4}, \quad\gamma_{3\cdot1}%
\]
one can construct many
other scalar differential invariants of PMAs. It can be shown these 
invariant does not
exhaust all invariants. Nevertheless, various quintuples of
(functionally) independent invariants can be chosen among them.
This is the only one need in order to apply the ``principle of
$n$ invariants''. For instance, we have

\begin{theorem}
Differential invariants composing each of the following two quintuples $\left(
\kappa_{1},\kappa_{2},\tau_{1},\tau_{2},\gamma_{3}\right) $ and
$\left(  \gamma_{3},N_{1}(\gamma_{3}),N_{1}^{2}(\gamma_{3}),\kappa
_{1},\kappa_{2}\right)  $ are independent.
\end{theorem}

\begin{proof}
 It is sufficient to exhibit an example for which
invariants composing each of these two quintuples are independent.
For instance,  the PMA with
directing distribution $\mathcal{R}$
generated by the field%
\[
Z=q\partial_{x}+y\partial_{y}+(qp+yq)\partial_{u}+x\partial_{p}-xu\partial
_{q}.
\]
is a such one.
In this case the corresponding Lagrangian distribution
$\mathcal{D}$ is generated by
\[
X_{1}=a_{1}(\partial_{x}+p\partial_{u})+a_{2}\partial_{p}+a_{3}\partial
_{q},\qquad X_{2}=a_{1}(\partial_{y}+q\partial_{u})+a_{3}\partial_{p}%
+a_{4}\partial_{q}%
\]
with
\[%
\begin{array}
[c]{l}%
a_{1}=xyu-(x-y)q,\qquad a_{2}=-(x-y)x-(u+xp)y^{2},\\
a_{3}=x^{2}u+(u+xp)yq,\qquad a_{4}=-(q+xu)xu-(u+xp)q^{2},
\end{array}
\]
and the corresponding PMA is%
\begin{equation}%
\begin{array}
[c]{l}%
a_{1}^{2}(rt-s^{2})+a_{1}a_{4}r+2a_{1}a_{3}s-a_{1}a_{2}t\\
\qquad+xa_{1}(xyup-xu-xpq+yu^{2}-uq)=0.
\end{array}
\label{example}%
\end{equation}

Explicit expressions of invariants $\left(
\kappa_{1},\kappa_{2},\tau _{1},\tau_{2},\gamma_{3}\right)  $ and
$\left(\gamma_{3},N_{1}(\gamma _{3}),\right. \allowbreak\left.
N_{1}^{2}(\gamma_{3}), \kappa_{1},\kappa _{2}\right)$ for equation
(\ref{example}) are too cumbersome to be reported here. At this point independence of invariants composing each 
of the above two quintuples is proved by a direct computer supported computation.
\end{proof}

\begin{remark}\label{complessity}
\begin{enumerate}
\item A proof as above based on computer computations could create a feeling of dissatisfaction. 
However, the use of computers for theoretical purposes is unavoidable in many questions,
say,  symmetries, conservation laws, recursion operators, etc, related with "sufficiently" 
nonlinear PDEs (see \cite{MVY,VK,KK}). A PCB composed of mutually projectively 
nonequivalent curves gives an informal geometrical explanation of this fact 
in the current context.
\item By consulting \cite{MVY,DPV} the reader will get a more precise idea on 
complexity of scalar differential invariants of MA equations.
\end{enumerate}
\end{remark}

Thus, according to the \textquotedblleft principle of $n$%
-invariants\textquotedblright\ (see \cite{AVL,SDI}), the proven
existence of five independent scalar differential invariants
solves in principle the equivalence problem for generic PMA
equations. It should be, however, stressed that a practical
implementation of this result could meet some boring/``unhuman" 
computations which, in fact, reflect the complexity of truly nonlinear PMAs.\\

\section{Examples}

Two examples in this section, first, illustrate complexity of coordinate expressions of invariants introduced in the previous section, and, second, show that these invariants are not necessarily independent. Also , the reader may observe that coordinate expressions of directing distributions are much simpler than those of corresponding generic PMAs. This gives an experimental evidence to the idea that working with directing distributions rather than original PMAs is more convenient.
Finally, it is worth noticing that directing distributions in these examples are
very simple, while their coordinate expressions are rather complicated.
\begin{example} 
The PMA
$\vspace{2mm}$

\[
\begin{array}{l}
\left(-xu+x^{2}p+xuq\right)\left(rt-s^{2}\right)+\left(x-qx\right)r\vspace{5pt}\\
\qquad\qquad+2xps+\left(-xp-uq+pu+u\right)t-p+q-1=0\end{array}\]

$\vspace{5mm}$

is generic. Its directing distribution is generated by
\[
Z=x(\partial_{x}+p\partial_{u})+u(\partial_{y}+q\partial_{u})+\partial_{p}+\partial_{q}.\]

The invariants of this equation are functions of $q$ and hence are not independent.
\[
\begin{array}{l}
k_{1}=\left[-3\,\left(-2-q\right)^{4}-{\frac{117}{5}}+4\, q\left(-2-q\right)+12\, q\left(-2-q\right)^{2}\right.\vspace{5pt}\\
\qquad\left.-24\, q+12\,\left(-2-q\right)^{2}-{\frac{81}{25}}\,\left(2\, q-\left(-2-q\right)^{2}\right)^{2}\right]\times\vspace{5pt}\\
\qquad\times\left(2\, q-6\, q\left(-2-q\right)+2\,\left(-2-q\right)^{3}+3\right)^{-4/3},\end{array}\]

$\vspace{5mm}$

\[
\begin{array}{l}
k_{2}=\left[\left(12\, q-36\, q\left(-2-q\right)+12\,\left(-2-q\right)^{3}+18\right)\left(-12-12\, q\right)\right.\vspace{5pt}\\
\qquad-7\,\left(14+12\, q-6\,\left(-2-q\right)^{2}\right)^{2}-\left({\frac{216}{5}}\, q-{\frac{108}{5}}\,\left(-2-q\right)^{2}\right)\times\vspace{5pt}\\
\qquad\left.\times\left(2\, q-6\, q\left(-2-q\right)+2\,\left(-2-q\right)^{3}+3\right)^{2}\right]\times\vspace{5pt}\\
\qquad\times\left(2\, q-6\, q\left(-2-q\right)+2\,\left(-2-q\right)^{3}+3\right)^{-8/3},\end{array}\]

$\vspace{5mm}$

\[
\begin{array}{l}
\tau_{1}=\left[-3\,\left(q-2\right)^{4}+6\,\left(-q+3\right)\left(q-2\right)^{2}-12\, q+{\frac{138}{5}}\right.\vspace{5pt}\\
\qquad\left.-12\,\left(q-2\right)^{2}-{\frac{81}{25}}\,\left(-q+3-\left(q-2\right)^{2}\right)^{2}\right]\times\vspace{5pt}\\
\qquad\times\left(-3\,\left(q-2\right)\left(-q+3\right)+2\,\left(q-2\right)^{3}-9/2+3\, q\right)^{-4/3}.\end{array}\]

$\vspace{5mm}$

\[
\begin{array}{l}
\tau_{2}=\left[\left(-18\,\left(q-2\right)\left(-q+3\right)+12\,\left(q-2\right)^{3}-27+18\, q\right)\left(-18+12\, q\right)\right.\vspace{5pt}\\
\qquad-7\,\left(6\, q-12+6\,\left(q-2\right)^{2}\right)^{2}-\left(-{\frac{108}{5}}\, q+{\frac{324}{5}}-{\frac{108}{5}}\,\left(q-2\right)^{2}\right)\times\vspace{5pt}\\
\qquad\left.\left(-3\,\left(q-2\right)\left(-q+3\right)+2\,\left(q-2\right)^{3}-9/2+3\, q\right)^{2}\right]\times\vspace{5pt}\\
\qquad\times\left(-3\,\left(q-2\right)\left(-q+3\right)+2\,\left(q-2\right)^{3}-9/2+3\, q\right)^{-8/3},\end{array}\]

$\vspace{5mm}$

\[
\gamma_{3}=\frac{2\, q-6\, q\left(-2-q\right)+2\,\left(-2-q\right)^{3}+3}{-3\,\left(q-2\right)\left(-q+3\right)+2\,\left(q-2\right)^{3}-9/2+3\, q}.\]

$\vspace{5mm}$

The invariant vector field $N_{1}$ is:\[
N_{1}=\frac{1}{\sqrt[3]{-3\,\left(q-2\right)\left(-q+3\right)+2\,\left(q-2\right)^{3}-9/2+3\, q}}Z.\]
\end{example}

\vspace{3mm}

\begin{example}
The PMA
$\vspace{2mm}$\[
\begin{array}{l}
\left(x^{2}-y^{2}\right)\left(rt-s^{2}\right)+\left(u(y-x)+y^{2}(p-q)\right)r\vspace{5pt}\\
+(2u(x-y)+2xy(p-q))s+\left(u(y-x)+x^{2}(p-q)\right)t\vspace{5pt}\\
+u(x+y)(q-p)=0\end{array}\]

$\vspace{5mm}$

is generic. Its directing distribution is generated by
\[
Z=y(\partial_{x}+p\partial_{u})+x(\partial_{y}+q\partial_{u})+u\partial_{p}+u\partial_{q}.\]

The invariants $k_{1},k_{2},\tau_{1},\tau_{2}$ are all functions of $\gamma_3$. Namely, 
\[
\begin{array}{l}
k_{1}=-\frac{{210}^{2/3}}{1102500}\,\frac{\left(-26490\gamma_{3}+1316225\,\gamma_{3}^{2}+729\right)}{\left(3+25\,\gamma_{3}\right)^{2/3}\gamma_{3}^{4/3}},\vspace{5pt}\end{array}\]

$\vspace{5mm}$

\[
\begin{array}{l}
k_{2}=-{\frac{\sqrt[3]{210}}{1050}}\,\frac{\left(-309+12545{\gamma}_{3}\right)}{\gamma_{3}^{2/3}(3+25\gamma_{3})^{1/3}},\end{array}\]

$\vspace{5mm}$

\[
\begin{array}{l}
\tau_{1}={\frac{1}{157500}}\,{\frac{\left(-22887-119730\,\gamma_{3}+364825\,\gamma_{3}^{2}\right){210}^{2/3}}{\left(3+25\,\gamma_{3}\right)^{2/3}},}\vspace{5pt}\end{array}\]

$\vspace{5mm}$

\[
\begin{array}{l}
\tau_{2}=-{\frac{1}{9450}}\,{\frac{\left(-124470\,\gamma_{3}+121875\,\gamma_{3}^{2}+546875\,\gamma_{3}^{3}-32076\right)\sqrt[3]{210}}{\left(3+25\,\gamma_{3}\right)^{1/3}},}\vspace{5pt}\end{array}\]

$\vspace{5mm}$

where\[
\gamma_{3}={\frac{-3/2\, x-3/2\, y}{{\frac{25}{2}}\, x+70+{\frac{25}{2}}\, y}}.\]

$\vspace{5mm}$

The invariant vector field $N_{1}$  for this equation is\[
N_{1}=\frac{1}{\sqrt[3]{{\frac{25}{2}}\, x+70+{\frac{25}{2}}\, y}}Z.\]
\end{example}

We conclude this section by illustrating the difference between generic and special PMAs.
\begin{example} 
Lagrangian distributions for the heat equation $u_{xx}-u_{y}=0$ and the Burgers equation  $u_{xx}+uu_{x}-u_{y}=0$ are 

\[
\mathcal{D}=\left\langle \partial_{x}+p\partial_{u}+q\partial_{p}, \partial_{q}\right\rangle, \]

and 
\[
\mathcal{D}=\left\langle \partial_{x}+p\partial_{u}+(q-up)\partial_{p}, \partial_{q}\right\rangle,\] 

respectively, while $\left\langle\partial_q\right\rangle$ is the generating distribution for each of them.
\end{example}

\section{Concluding remarks}

Representation of a PMA $\mathcal{E}$ by means of the
associated PCB makes clearly visible the nature of its
nonlinearities. For example, if all the curves of this bundle are
projectively nonequivalent to each other, then $\mathcal{E}$ does not
admit contact symmetries. This occurrence can be detected by means of invariants constructed in section 7. 
On the other hand, it may happen that all curves composing
a PCB are projectively equivalent, i.e., nonlinearities of the
corresponding PMA $\mathcal{E}$ are ``homogeneous''. The invariants considered in this paper are not sufficient to distinguish two PCBs which are homogeneous in this sense. So, the need of some 
new invariants arises. It is remarkable that PCBs immediately suggest
an idea of how the necessary new invariants
can be constructed. For instance, one can observe that in this case 
the bundle $PT^*N\rightarrow N$ is naturally
supplied with a full parallelism structure which immediately
furnishes the required new invariants. It is not difficult to
imagine various intermediate situations, which demonstrate the
diversity and complexity of the world of PMAs. In particular, the problem 
of describing all strata of the characteristic diffiety (see \cite{SDI}) for 
PMAs appears to be a task of a rather large scale. Further
results in this direction will appear in forthcoming
publications.

\ack The first author DCF acknowledge the financial support of FAPESB (Funda\c c\~ao de Amparo \`a Pesquisa do Estado da Bahia), which made possible the conclusion of this research. However, since a preliminary research was already achieved during his stay at Universit\`a Statale di Milano and Silesian University of Opava, DCF also acknowledge these two institutions and in particular Michal Marvan for his kind hospitality and financial support through grant MSM 4781305904.


\end{document}